\documentclass[10pt, amsmath, amsthm, amssymb]{amsart}

\usepackage{mathtools}

\usepackage{pdfpages}
\usepackage[left=20mm,top=0.6in,bottom=12mm]{geometry}

\allowdisplaybreaks[1]

\numberwithin{equation}{section}

\newtheorem{Lemma}{LEMMA}[section]
\newtheorem{Theorem}[Lemma]{Theorem}
\newtheorem{Proposition}[Lemma]{Proposition}
\newtheorem{Corollary}[Lemma]{Corollary}

\newtheorem{remark}[Lemma]{Remark}
\newtheorem{definition}[Lemma]{Definition}
\newtheorem{example}[Lemma]{Example}

\newtheorem{Fact}[Lemma]{Fact}
\newtheorem{assumption}[Lemma]{Assumption}

\def\bt{\begin{Theorem}}
\def\et{\end{Theorem}}
\def\bl{\begin{Lemma}}
\def\el{\end{Lemma}}
\def\bp{\begin{Proposition}}
\def\ep{\end{Proposition}}
\def\bcor{\begin{Corollary}}
\def\ecor{\end{Corollary}}
\def\bpf{\begin{proof}}
\def\epf{\end{proof}}

\def\brem{\begin{remark}\rm }
\def\erem{\hfill $\lozenge$ \end{remark}}

\def\bedef{\begin{definition}\rm }
\def\endef{\hfill$\lozenge$\end{definition}}

\def\beg{\begin{example}\rm }
\def\eeg{\hfill $\lozenge$\end{example}}

\def\bef{\begin{Fact}}
\def\eef{\end{Fact}}

\def\bea{\begin{assumption}}
\def\ena{\end{assumption}}
\def\bc{\begin{center}}
\def\ec{\end{center}}
\def\noi{\noindent}

\def\vsq{\vskip .25cm}
\def\beq{\begin{equation}}
\def\eeq{\end{equation}}
\def\beqarray{\begin{eqnarray*}}
\def\eeqarray{\end{eqnarray*}}
\def\<{\leftangle}
\def\>{\rightangle}
\def\({\left(}
\def\){\right)}
\def\f{\varphi}

\def\<{\langle}
\def\>{\rangle}
\def\q{\quad}
\def\r{\rho}
\def\a{\alpha}

\def\d{\delta}

\def\k{\kappa}

\def\t{\tau}
\def\r{\rho}
\def\l{\lambda}
\def\e{\varepsilon}
\def\s{\sigma}
\def\O{\Omega}
\def\o{\omega}

\def\m{\mu} 

\def\w.r.t.{with respect to}
\def\R{{\mathbb{R}}}
\def\N{{\mathbb{N}}}

\def\bq{\begin{quote}}
\def\eq{\end{quote}}

\def\bit{\begin{itemize}}
\def\eit{\end{itemize}}
\def\i{\item}
\def\ben{\begin{enumerate}}
\def\een{\end{enumerate}}

\def\ds{\displaystyle}

\begin{document}
\title[inverse source identification for bi-parabolic equation]{A source identification problem in a bi-parabolic equation: convergence rates and some optimal results}
\author{Subhankar Mondal and M. Thamban Nair}
\address{(SM): TIFR Centre for Applicable Mathematics, Bangalore-560065, India;
(MTN): Department of Mathematics, BITS Pilani, Goa campus, Goa-403726}
%Institute for Numerical and Applied Mathematics, University of G{\"o}ttingen, G{\"o}ttingen-37083, Germany}
\email{ subhankar22@tifrbng.res.in; 
mtnair@goa.bits-pilani.ac.in
}
\maketitle

\begin{abstract} 
This paper is concerned with identification of a spatial source function from final time observation in a bi-parabolic equation, where the full source function is assumed to be a product of time dependent and a space dependent function. Due to the ill-posedness of the problem, recently some authors have employed different regularization method and analysed the convergence rates. But, to the best of our knowledge, the quasi-reversibility method is not explored yet, and thus we study that in this paper. As an important implication, the H{\"o}lder rates for the apriori and aposteriori error estimates obtained in this paper improve upon the rates obtained in earlier works. Also, in some cases we show that the rates obtained are of optimal order. Further, this work seems to be the first one that has broaden the applicability of the problem by allowing the time dependent component of the source function to change sign. To the best of our knowledge, the earlier known work assumed the fixed sign of the time dependent component by assuming some bounded below condition.
\end{abstract}


\section{Introduction} 
In the theory of heat transfers, the bi-parabolic equations have many applications \cite{payne_song_2006}. It is reported, for instance in \cite{fichera_1992, joseph_preziosi_1989}, that the classical parabolic equations fail to  describe the heat conduction processes accurately. To overcome this shortcoming, several models were proposed which include the bi-parabolic model proposed in \cite{fushchich_galitsyn_polubinskii_1990}. In addition to this, bi-parabolic models are used to describe the special phenomenon of process dynamics filter image \cite{kalantarov_zelik_2009} and also to model various natural phenomena which includes ice formation, fluid on the lungs, design of distinctive curves on the surface \cite{greer_bertozzi_sapiro_2006}. Thus, from application point of view the bi-parabolic systems are interesting and relevant topic.

In this paper, we consider an inverse source identification problem associated with a bi-parabolic system. More precisely, let $\O\subset \R^d (d\geq 1),$ be a bounded domain with Lipschitz boundary $\partial\O$ and $\t>0$ be fixed. We consider the following system
\beq\label{biparbolic_exact}
\begin{cases}
u_{tt}-2\Delta u_t +\Delta^2u=\psi(t)f(x),&\q\text{in}\q\O\times (0,\t),\\
u=0=\Delta u,&\q\text{on}\q\partial\O\times (0,\t),\\
u(\cdot,0)=0=u_t(\cdot,0),&\q\text{in}\q\O,\\
u(\cdot,\t)=h,&\q\text{in}\q\O,
\end{cases}
\eeq
where $\psi\in C([0,\t])$ and $f, h\in L^2(\O).$ The inverse problem is to identify the spatially dependent exact source $f$ from the knowledge of exact final time observation $h.$ But in practice since we will always have access to a measurement data of the final value, noise in the data is inevitable. In this paper, we consider the deterministic noise model. That is, for $\d>0,$ let $h^\d$ be the measured final time noisy observation satisfying 
\beq\label{noise_model}
\|h-h^\d\|_{L^2(\O)}\leq \d.
\eeq
This type of problem belongs to the class of inverse source identification problem which has a vast literature in the area of inverse problems (cf. \cite{isakov_book}) because of their tremendous real world applications. Indeed, the source identification problem is used in detecting the source of pollution in water surfaces or atmospheric media \cite{andrle_elbadia_2012}, and also for geophysical prospecting \cite{magnoli_viano_1997}, to name a few.

The inverse source problem that is considered in this paper is known to be ill-posed \cite{quoc_long_oregan_ngoc_tuan_2022}, in the sense that a small perturbation in the final time observation may lead to a large deviation in the corresponding source function and thus, some regularization method (cf. \cite{engl_hanke_neubauer, nair_opeq}) has to be employed in order to obtain some stable approximations for the source function. For $\psi\equiv 1$ in $[0,\t]$, an iterative regularization method has been proposed in \cite{zouyed_djemoui_2015} to obtain stable approximations for the identification of the spatially dependent source function. For the general source term of the separable form that is considered in this paper, it seems the work on the regularization aspect is in the initial stage. To the best of our knowledge, the very recent paper of Quoc Nam et al. \cite{quoc_long_oregan_ngoc_tuan_2022} is the first paper where the authors considered the same problem of source identification that is considered here, and  have employed the Tikhonov regularization method for the regularization purpose. Subsequently, in a follow-up work, Tuan \cite{tuan_2022} has used a Fourier truncation method for the regularization purpose again for a similar source identification problem and analysed the case when the observed data is in any $L^p$ spaces. In \cite{phuong_luc_long_2020} the authors employed the quasi-boundary value method for the regularization purpose for the source identification problem that is considered in this paper. To the best of our knowledge, another important type of regularization method, namely the quasi reversibility method (QRM) has not been explored yet in the setting of source identification associated with the bi-parabolic equation \eqref{biparbolic_exact}. In this paper we employ the QRM for the regularization purpose and obtain the error estimates.

The QRM was popularised by Latt{\`e}s and Lions \cite{lattes_lions_1967}. The idea behind this method is to perturb the governing differential equation by a suitable operator so that solving the resulting problem becomes a well-posed problem. More specifically, in \cite{lattes_lions_1967} the authors considered a final value problem of the form
$$
\begin{cases}
\partial_t v+ \mathcal{L}(t)v =0\\
 v(\t)= v_\t
\end{cases}
$$
and considered the problem of identifying the initial value $v(0),$ where the operator $\mathcal{L}(t)$ is a spatial differential operator, for example, one may have $\mathcal{L}(t)$ to be the negative spatial Laplacian $-\Delta$ or of the form $-\nabla\cdot a(x,t)\nabla$, where $a(x,t)$ is sufficiently smooth. In order to solve the backward problem the authors proposed the following perturbed equation
$$
\begin{cases}
\partial_t v_\epsilon+\mathcal{L}(t)v_\epsilon+\epsilon \mathcal{L}^*(t)\mathcal{L}(t)v_\epsilon=0\\
 v_\epsilon(\t)= v^\epsilon_\t
\end{cases}$$
for $\epsilon>0$. The problem of obtaining the initial value $v_\epsilon(0)$ from the final value $v^\epsilon_\t$ is a well-posed problem and $v_\epsilon(0)$ gives an approximation of $v(0)$ in some suitable topology. Since the work in \cite{lattes_lions_1967}, several researchers have analysed the QRM in a wide range of problems. The application of QRM can be found, for example, in \cite{bourgeois_2006} for the elliptic setting, \cite{dorroh_ru_1999, mophou_warma_2023, nguyen_khoa_vo_2019, showalter_1975, tuan_nane_trong_2021} for the parabolic setting and \cite{le_nguyen_nguyen_powell_2021, nguyen_2019} for the hyperbolic setting. Also, with the increasing work in the backward time fractional diffusion equation, the QRM has found its application in several recent works related to time/space fractional setting, see for instance \cite{duc_thang_thanh_2023, yang_ren_li_2018}. For more details on QRM and its several implications and modifications, the interested reader may refer to  \cite{ames_clark_epperson_oppenheimer_1998, clark_oppenheimer_1994, clason_klibanov_2007, huang_2008, showalter_1974}.

Motivated by the recent work of Duc et al. \cite{duc_thang_thanh_2023}, where a QRM is analysed for a source identification problem associated with a time-space fractional parabolic equation, we analyse a QRM for the inverse problem considered in this paper. More precisely, for $b\geq 2$ and $\a>0,$ we consider the following bi-parabolic system
$$
\begin{cases}
v_{tt}-2\Delta v_t +\Delta^2v=(I+\a(-\Delta)^b)\psi(t)\tilde{f}(x),&\q\text{in}\q\O\times (0,\t),\\
v=0=\Delta v,&\q\text{on}\q\partial\O\times (0,\t),\\
v(\cdot,0)=0=v_t(\cdot,0),&\q\text{in}\q\O,\\
v(\cdot,\t)=\tilde{h},&\q\text{in}\q\O,
\end{cases}
$$
where $I:L^2(0,\t;L^2(\O))\to L^2(0,\t;L^2(\O))$ is the identity map, and consider the problem of identifying $\tilde{f}$ from the final time measurement $\tilde{h}$ using the {\it mild solution} (see Definition \ref{mild_sol_qrm}). We show that this is an well-posed problem, and the obtained solutions are indeed an approximations for the exact source $f.$

We now point out the novelty and main contributions of this paper.
\bit
\i In all the earlier known works \cite{quoc_long_oregan_ngoc_tuan_2022, tuan_2022, phuong_luc_long_2020} that considered the identification of $f$ from $h$ associated with the bi-parabolic system \eqref{biparbolic_exact}, the analysis has been done by assuming a bounded below assumption on $\psi,$ the time dependent component of the source, which ensures a fixed sign property for $\psi$. In contrast to this, our analysis works for a more general case, which includes the fixed sign case and also allows $\psi$ to change sign (see Assumption \ref{assumption_main}). In that sense, this work seems to be the first paper that allows $\psi$ to change sign in the context of the considered inverse problem.
\i This work seems to be the first work where the QRM is analysed for the regularization purpose in the context of the bi-parabolic system \eqref{biparbolic_exact}. As mentioned earlier, the Tikhonov regularization method, quasi-boundary value method and the Fourier truncation method have already been analysed by some authors for the same problem.
\i We obtain H{\"o}lder type error estimates for both the apriori and aposteriori parameter choice strategies using a conditional stability estimates. Under some additional assumptions on $b$, the order of the error estimates in apriori case is better than the previously known rates (see Remarks \ref{better_rate_apriori} and \ref{better_rate_apost}).
\i We show that in some cases the obtained error estimates are of optimal order (in the sense of {\it worst case error}) for both the apriori and aposteriori cases. This fills the gap in the arguments for the order optimality claim made in \cite{duc_thang_thanh_2023, hao_liu_duc_thang_2019} (see Remark \ref{incomp_arg_duc_hao}).
\eit
The remainder of the paper is organized as follows. In Section \ref{sec_prelim} we recall all the necessary results and state some facts that will be used throughout in this paper. Section \ref{sec-assum_cond_est} perhaps is the base of the paper. In this section, we specify certain assumptions on $\psi$, see their consequences and as a by product we obtain a conditional stability estimate. In Section \ref{sec-reg} we consider the regularization by quasi-reversibility method, obtain error estimates for both apriori and aposteriori parameter choice strategies. Finally, in Section \ref{sec-optimality} we prove that some of the error estimate that we obtain are of optimal order.

\section{preliminaries}\label{sec_prelim}
In this section we recall some basic results from the spectral theory of elliptic PDE, and discuss about the solutions associated to \eqref{biparbolic_exact}, and also present some basic results that will be used later.

First, let us fix some notations that will be used throughout.
\bit
\i Let $\mathbb{X}$ denote the space of real-valued functions defined on $\O$, for e.g., $L^2(\O).$ Then, for $z\in L^2(0,\t;\mathbb{X}),$ we may use the notation $z(t)$ to denote $ z(\cdot,t)$ for almost all $t\in [0,\t].$ %\hfill {\bf to be re-written} 
\i $C([0,\t];L^2(\O))$ denotes the standard space of all $L^2(\O)$-valued continuous function in $[0,\t].$
\i If $\chi\in C([0,\t])$, then $\|\chi\|_\infty$ denotes the standard $\sup$ norm.
\i For $g_1,g_2\in L^2(\O)$, $\<g_1,g_2\>$ denotes the standard inner-product in $L^2(\O).$
\i For normed linear spaces $\mathbb{X}$ and $\mathbb{Y}$, and a linear operator $\mathbb{A}:\mathbb{X}\to\mathbb{Y},$ $\text{\bf dom}(\mathbb{A})$, $\text{\bf ran}(\mathbb{A})$ and $N(\mathbb{A})$ shall denote the domain, range and the kernel of $\mathbb{A},$ respectively.
\eit
For $\lambda\in \R,$ it is known that the spectral problem
\beq\label{eigen_value_Laplacian}
\begin{cases}
-\Delta v=\lambda v\q &\text{in}\q\O,\\
\q v=0\q &\text{on}\q\partial\O,
\end{cases}
\eeq
 admits a sequence of eigenvalues satisfying $0<\lambda_1\leq \lambda_2\leq\ldots\to +\infty$, and the corresponding eigenfunctions $\f_n$ are such that $\f_n\in H^1_0(\O)\cap H^2(\O)$ and $\{\f_n:n\in \N\}$ is an orthonormal basis of $L^2(\O)$ (see for e.g., \cite{evans}).
 
For $p\geq 0,$ we define 
$$\mathbb{H}_p:=\{\phi\in L^2(\O):\sum_{n=1}^\infty \lambda_n^{2p}|\<\phi,\f_n\>|^2<\infty\},$$ 
where $\<\cdot, \cdot\>$ denotes the inner product  $\<\cdot,\cdot \>_{L^2(\O)}$ on the space $L^2(\O)$.  
It can be seen that  $\mathbb{H}_p$ is a Hilbert space \w.r.t. the inner product
$$\<\phi, \psi\>_{{\mathbb H}_p}  = \sum_{n=1}^\infty \lambda_n^{2p}\<\phi,\f_n\> \<\f_n, \psi\>,\q  \phi, \psi\in {\mathbb H}_p,$$
so that the corresponding   norm is given by  $$\|\phi\|_{\mathbb{H}_p}=\left(\sum_{n=1}^\infty \lambda_n^{2p}|\<\phi,\f_n\>|^2\right)^\frac{1}{2}, \q  \phi \in {\mathbb H}_p.$$
Clearly, when $p=0,$ then $\mathbb{H}_p=L^2(\O).$

We now look into the existence and uniqueness of the solution of the following initial boundary value problem which will help to obtain a representation of the source term $f$ in terms of $\psi$ and $h$. Let $\psi\in C([0,\t])$, $f\in L^2(\O)$. We consider the following system
\beq\label{biparbolic_forward}
\begin{cases}
v_{tt}-2\Delta v_t +\Delta^2v=\psi(t)f(x),&\q\text{in}\q\O\times (0,\t),\\
v=0=\Delta v,&\q\text{on}\q\partial\O\times (0,\t),\\
v(\cdot,0)=0=v_t(\cdot,0),&\q\text{in}\q\O.
\end{cases}
\eeq
%Following \cite{zouyed_djemoui_2015, quoc_long_oregan_ngoc_tuan_2022} (see also \cite{tuan_2022}) the problem \eqref{biparbolic_forward} has a mild solution given by
%\beq\label{sol_rep_exact_biparabolic}
%u(t)=\sum_{n=1}^\infty u_n(t)\f_n,
%\eeq
%where {\cb $u_n(t):=\<u(t),\f_n\>$ is given by}
%$$u_n(t)=(1+\lambda_nt)e^{-\lambda_nt}u_n(0)+{\cb f_n\int_0^t e^{-\lambda_n(t-s)}(t-s)\psi(s)\,ds},$$
%and $\,f_n=\<f,\f_n\>.$ 
%Following \cite{zouyed_djemoui_2015, quoc_long_oregan_ngoc_tuan_2022} (see also \cite{tuan_2022}) the problem \eqref{biparbolic_forward} has a solution {\cb $u,$} given by
%\beq\label{sol_rep_exact_biparabolic}
%u(t)=\sum_{n=1}^\infty u_n(t)\f_n,
%\eeq
%with 
%$$u_n(t):= \<u(t),\f_n\> =  \<f, \f_n\> \int_0^t e^{-\lambda_n(t-s)}(t-s)\psi(s)\,ds.$$
Following \cite{zouyed_djemoui_2015, quoc_long_oregan_ngoc_tuan_2022} (see also \cite{tuan_2022}) the problem \eqref{biparbolic_forward} has a solution $u,$ given by
\beq\label{sol_rep_exact_biparabolic}
u(t)=\sum_{n=1}^\infty \mu_n(t) \<f, \f_n\> \f_n,
\eeq
with 
\beq\label{ev-t}\mu_n(t):=  \int_0^t e^{-\lambda_n(t-s)}(t-s)\psi(s)\,ds.\eeq

%\bedef\label{mild_sol_exact}\cb
%For $\psi\in L^2(0,\t)$ and $f\in L^2(\O)$, a function $u\in C([0,\t];L^2(\O))$ with the representation as given in \eqref{sol_rep_exact_biparabolic} is said to be a mild solution of \eqref{biparbolic_forward}.
%\endef
%\brem\cb
%For $\psi\in L^2(0,\t)$ and $f\in L^2(\O),$ the mild solution $u$ of \eqref{biparbolic_forward} is also a solution, that is, $u$ satisfies \eqref{biparbolic_forward} in the $L^2$ sense.
%\erem
Thus, taking into account the final value condition $u(\t)=h$ in \eqref{biparbolic_exact}, we have 
\beq\label{rep-h}\<h,\f_n\>=\int_0^\t e^{-\lambda_n(\t-s)}(\t-s)\psi(s)\,ds\,\<f,\f_n\> = \m_n \<f, \f_n\>,\eeq 
where 
\beq\label{ev}  \m_n:= \m_n(\t) = \int_0^\tau e^{-\lambda_n(\t-s)}(\t-s)\psi(s)\,ds.\eeq 
In particular, the inverse source identification problem associated with \eqref{biparbolic_exact}  is same as the problem of solving the operator equation 
\beq\label{op-eq} {\mathbb T}f = h,\eeq 
where ${\mathbb T} : L^2(\O)\to L^2(\O)$ is defined by 
\beq\label{cpt-op} {\mathbb T}\f = \sum_{n=1}^\infty \mu_n \<\f, \f_n\>\f_n,\q \f\in L^2(\O).\eeq 
Since $\l_n\to \infty$ and $\psi\in C[0, \t]$, it can be seen that  $|\mu_n|\to 0$ as $n\to \infty$.  In fact, 
	\beq\label{bound_for_int_e_psi}  
	|\mu_n|\leq \|\psi\|_{\infty} \int_0^\t e^{-\lambda_n(\t-s)}(\t-s)\,ds=\|\psi\|_{\infty}\frac{1-(1+\t\lambda_n)e^{-\lambda_n\t}}{\lambda_n^2}
	\leq  \|\psi\|_{\infty}\,\frac{1}{\lambda_n^2}.
	\eeq
Hence,  from the representation (\ref{cpt-op}) of $\mathbb T$, it follows that  ${\mathbb T}$ is a compact self-adjoint operator with eigenvalues $\mu_n,\, n\in \N.$ In particular, we can infer that the inverse problem under consideration, which is same as the problem of solving the operator equation \eqref{op-eq} is an ill-posed problem.
We shall show, under certain assumptions on $\psi$, that $\m_n\not= 0$ for all $n\in \N$ so that under those conditions, the equation \eqref{op-eq} has a unique solution provided 
$$\sum_{n=1}^\infty \frac{|\<h, \f_n\>|^2}{|\m_n|^2}<\infty,$$
and in that case the solution $f$ is given by 
\beq\label{rep_f_intermsof_h} f = \sum_{n=1}^\infty \frac{\<h, \f_n\>}{\m_n} \f_n. \eeq
%
%\beq\label{rep_f_intermsof_h}
%\<f,\f_n\>=\frac{\<h,\f_n\>}{\int_0^\t e^{\lambda_n(s-\t)}(\t-s)\psi(s)\,ds}.
%\eeq
%In order that \eqref{rep_f_intermsof_h} makes sense, we must ensure that the denominator is non zero. For this we make certain assumptions which will {\cb be} discussed in the next section.
%We now state some useful results that will be used in the later part.
%\bl\label{eta_a}
%Let $a>0$ be fixed and $\eta_a:[\lambda_1,\infty)\to \R$ be defined by $\eta_a(s)=(1+as)e^{-as}$. Then $\eta_a$ is strictly decreasing and hence attains its maximum at $s=\lambda_1.$
%\el
%\bpf
%The proof follows by observing $$\frac{d\eta_a}{ds}=ae^{-as}-a(1+as)e^{-as}=-a^2se^{-as}<0.$$
%\epf
%\bl\label{eta_alpha_b_p}
%Let $\a, p>0$, and $b>p.$ Let $\eta_{\a,b,p}:[0,\infty)\to [0,\infty)$ be defined by $\eta_{\a,b,p}(s)=\frac{\a s^{b-p}}{1+\a s^b}.$ Then $\eta_{\a,b,p}$ attains global maximum at $s_0=\big(\frac{b-p}{\a p}\big)^\frac{1}{b}$ and $\eta_{\a,b,p}(s_0)=\frac{p}{b}\big(\frac{b-p}{p}\big)^\frac{b-p}{b}\,\a^\frac{p}{b}.$
%\el
%\bpf
%The proof follows easily by observing that {\cb $\eta_{\a,b,p}$ is differentiable and $\eta_{\a,b,p}'(s)$ vanishes only at $s_0.$}
%\epf

We now state two elementary results in the form of a lemma, whose proofs follow easily from basic calculus. 

\bl\label{lemma-basic}   
\ben 
\i[(i)] Let $a>0$ be fixed and $\eta_a:[\lambda_1,\infty)\to \R$ be defined by $\eta_a(s)=(1+as)e^{-as}$. Then $\eta_a$ is strictly decreasing and hence attains its maximum at $s=\lambda_1.$
\i[(ii)] Let $\a, p>0$, and $b>p.$ Let $\eta_{\a,b,p}:[0,\infty)\to [0,\infty)$ be defined by $\eta_{\a,b,p}(s)=\frac{\a s^{b-p}}{1+\a s^b}.$ Then $\eta_{\a,b,p}$ attains global maximum at $s_0=\left(\frac{b-p}{\a p}\right)^\frac{1}{b}$ and $\eta_{\a,b,p}(s_0)=\frac{p}{b}\left(\frac{b-p}{p}\right)^\frac{b-p}{b}\,\a^\frac{p}{b}.$
\een \el

\section{assumption on $\psi$, its consequences and a stability result}\label{sec-assum_cond_est}
The results of this section are the base of this paper. We begin with the following assumption.
\bea\label{assumption_main}
There exists $\t_0\in [0,\t)$ and $\k_1>0$ such that $|\psi|\geq \k_1$ on $[\t_0,\t],$ and one of the following holds:
\ben
\i[(i)] $\psi$ does not change sign on $[0,\t].$
\i[(ii)] If $\psi$ changes sign on $[0,\t]$ then $\psi$ is differentiable and there exists $\k_2>0$ such that $|\psi'|\leq \k_2$ on $[0,\t],$ and $|\psi|\leq M:=\frac{ (\t-\t_0)^2\k_1}{4\t \t_0}$ on $I_{\text{sc}}:=\{t\in [0,\t]:\psi(t)\psi(\t)\leq 0\}.$
\een
\ena
\brem\label{remark_assump}
The above assumption on $\psi$ is less restrictive than the assumptions made in  \cite{quoc_long_oregan_ngoc_tuan_2022, tuan_2022, phuong_luc_long_2020} (where the analysis is done by assuming $|\psi|$ is bounded below by a positive constant, which by virtue of continuity implies $\psi$ has a fixed sign, i.e., either $\psi>0$ or $<0$ in $[0,\t]$) in the following sense: firstly, the bounded below assumption for the fixed sign case does not necessarily has to be on the whole interval $[0,\t]$, and secondly, the sign changing assumption obviously implies $\psi$ can take much general form than only the fixed sign case.  
\erem 
We now give an example of a function that satisfies Assumption \ref{assumption_main} (ii).
\beg
Let $\t=\pi.$ Consider the function $$\psi(t)=\begin{cases}
\frac{1}{48\sqrt{2}}\cos t,\q&\text{if}\q t\in[0,\frac{\pi}{2}],\\
\cos t,\q&\text{if}\q t\in [\frac{\pi}{2},\pi].
\end{cases}
$$ Clearly $|\psi(t)|\geq \frac{1}{\sqrt{2}}$ in $[\frac{3\pi}{4},\pi],$ and hence $\t_0:=\frac{3\pi}{4}$ and $\k_1:=\frac{1}{\sqrt{2}}.$ Since $\psi(\pi)<0$, $I_{sc}:=[0,\frac{\pi}{2}].$ Now, observe that $\max\limits_{I_{sc}}|\psi|=\frac{1}{48\sqrt{2}}\leq M:=\frac{ (\t-\t_0)^2\k_1}{4\t \t_0}=\frac{1}{48\sqrt{2}}.$ Thus, $\psi$ satisfies all the requirements of Assumption \ref{assumption_main} (ii).
\eeg

\brem
Although the Assumption \ref{assumption_main} is motivated from \cite[Assumption H]{duc_thang_thanh_2023}, there is a difference to that. The difference is w.r.t. the bound $M$ of $|\psi |$ on $I_\text{sc}.$ Note that $M= \frac{(\t-\t_0)}{4\t}M_0,$ where $M_0:=\frac{(\t-\t_0)\k_1}{\t_0}$ is the bound  in \cite{duc_thang_thanh_2023}.
\erem

\bt\label{non_zero_lower_bound}
Let Assumption \ref{assumption_main} holds and $\m_n$ be as in (\ref{ev}). Then there exists a constant $C>0$ such that 
$$
|\m_n|\geq \frac {C}{\l_n^2}\q\forall \, n\in \N.$$
In particular, $\m_n\not=0$ for all $n\in \N$. 
\et

\bpf
{\bf \underline{Case 1.}} Suppose the item (i) of Assumption \ref{assumption_main} holds. Then either $\psi(s)>0$ or $\psi(s)<0$ for all $s\in [0,\t].$ Since $e^{\lambda_n(s-\t)}(\t-s)\geq 0$ for $s\in [0,\t]$, we obtain 
$$|\m_n| = 
	\Big|\int_0^\t e^{\lambda_n(s-\t)}(\t-s)\psi(s)\,ds\Big|= \int_0^\t e^{\lambda_n(s-\t)}(\t-s)|\psi(s)|\,ds.$$ 
Note that 
\beqarray
\l_n^2 \int_0^\t e^{\lambda_n(s-\t)}(\t-s)|\psi(s)|ds &\geq &\k_1 \lambda_n^2\int_{\t_0}^\t e^{\lambda_n(s-\t)}(\t-s)\,ds\\
&=& \k_1[1-(1+\lambda_n(\t-\t_0))e^{-\lambda_n(\t-\t_0)}]\\
& \geq & \k_1[1-(1+\lambda_1(\t-\t_0))e^{-\lambda_1(\t-\t_0)}].
\eeqarray
The last inequality follows from the fact that  the function  $t\mapsto (1+t)e^{-t}$ is strictly  decreasing on $(0, \infty)$, thanks to Lemma \ref{lemma-basic}(i). Since $(1+\lambda_1(\t-\t_0))e^{-\lambda_1(\t-\t_0)} <1$, we obtain the required inequality with $C =  \k_1[1-(1+\lambda_1(\t-\t_0))e^{-\lambda_1(\t-\t_0)}]$. 
\vsq

\noi 
{\bf \underline{Case 2.}} Suppose that item (ii) of Assumption \ref{assumption_main} holds, and let  $\t_1=\max\{t\in [0,\t]:\psi(t)=0\}.$ It is easy to observe that $0<\t_1<\t_0.$ Let $I_{\text{su}}:=\{t\in [0,\t]:\psi(t)\psi(\t)\geq 0\}.$ Then $[\t_1,\t]\subset I_{\text{su}}$ and $I_{\text{sc}}\subset [0,\t_1].$ Let $\tilde{C}$ be such that $\tilde{C}>\frac{|\psi(0)|+2\t \k_2}{\t|\psi(\t)|}.$

\noi {\bf \underline{Sub-case 2 (i).}}\q $\lambda_n\leq \tilde{C}.$ 

\noi
Recall that by Assumption \ref{assumption_main} (ii) we have $|\psi|\leq M :=\frac{ (\t-\t_0)^2\k_1}{4\t\t_0}$ on $I_{\text{sc}}:=\{t\in [0,\t]:\psi(t)\psi(\t)\leq 0\}.$ Hence, 
\beqarray
\lambda_n^2|\m_n| &\geq & \lambda_n^2\Big[\int_{I_{\text{su}}}e^{-\lambda_n(\t-s)}(\t-s)|\psi(s)|\,ds-\int_{I_\text{sc}}e^{-\lambda_n(\t-s)}(\t-s)|\psi(s)|\,ds\Big]\\
&\geq &  \lambda_n^2\Big[\int_{\t_1}^\t e^{-\lambda_n(\t-s)}(\t-s)|\psi(s)|\,ds-M\int_0^{\t_1}e^{-\lambda_n(\t-s)}(\t-s)\,ds\Big]\\
&= & \lambda_n^2\Big[\int_{\t_1}^{\t_0} e^{-\lambda_n(\t-s)}(\t-s)|\psi(s)|\,ds +\int_{\t_0}^{\t} e^{-\lambda_n(\t-s)}(\t-s)|\psi(s)|\,ds
 -M\int_0^{\t_1}e^{-\lambda_n(\t-s)}(\t-s)\,ds\Big]\\
&\geq & \lambda_n^2\Big[\int_{\t_1}^{\t_0} e^{-\lambda_n(\t-s)}(\t-s)|\psi(s)|\,ds+\k_1\int_{\t_0}^\t e^{-\lambda_n(\t - s)}(\t-s)\,ds -M\int_0^{\t_1}e^{-\lambda_n(\t-s)}(\t-s)\,ds\Big].
\eeqarray
We now obtain some bounds for the integrals $$I_1:=\int_{\t_1}^{\t_0} e^{-\lambda_n(\t-s)}(\t-s)|\psi(s)|\,ds,\q I_2:=\int_{\t_0}^\t e^{-\lambda_n(\t-s)}(\t-s)\,ds,\q I_3:=\int_0^{\t_1}e^{-\lambda_n(\t-s)}(\t-s)\,ds.$$
It is easy to see that 
$$I_1\geq (\t-\t_0)e^{-\lambda_n(\t-\t_1)}\int_{\t_1}^{\t_0}|\psi(s)|\,ds.$$
Let $\e=\frac{\t-\t_0}{2}$. Then
$$
I_2\geq \int_{\t_0}^{\t-\e}e^{-\lambda_n(\t-s)}{ (\t-s)}\,ds\geq \frac{\e{ (\t-\t_0-\e)}}{e^{\lambda_n(\t-\t_0)}}\geq \frac{\e{ (\t-\t_0-\e)}}{e^{\lambda_n(\t-\t_1)}} {= \frac{(\t-\t_0)^2}{4e^{\lambda_n(\t-\t_1)}}}.
$$
Similar to $I_1$ it is easy to check that $I_3\leq \t\t_0 e^{-\lambda_n(\t-\t_1)}.$ Thus, using these estimates we obtain
\beqarray
\lambda_n^2|\m_n| &\geq & \lambda_n^2[\,(\t-\t_0)e^{-\lambda_n(\t-\t_1)}\int_{\t_1}^{\t_0}|\psi(s)|\,ds+ { \frac{(\t-\t_0)^2\k_1}{4}} e^{-\lambda_n(\t-\t_1)}-M\t{ \t_0} e^{-\lambda_n(\t-\t_1)}\,]\\
&\geq &\lambda_n^2 e^{-\tilde{C}(\t-\t_1)} [\,(\t-\t_0)\int_{\t_1}^{\t_0}|\psi(s)|\,ds+{ \frac{(\t-\t_0)^2\k_1}{4}}-M\t{\t_0}\,]\\
&\geq & \lambda_1^2 e^{-\tilde{C}(\t-\t_1)} [\,(\t-\t_0)\int_{\t_1}^{\t_0}|\psi(s)|\,ds+{ \frac{(\t-\t_0)^2\k_1}{4}}-M\t{ \t_0}\,]\\
&= & \lambda_1^2 e^{-\tilde{C}(\t-\t_1)} (\t-\t_0)\int_{\t_1}^{\t_0}|\psi(s)|\,ds.
\eeqarray
Thus, the required inequality is satisfied with 
$C = \lambda_1^2 e^{-\tilde{C}(\t-\t_1)} (\t-\t_0)\int_{\t_1}^{\t_0}|\psi(s)|\,ds.$

\vsq
\noi 
{\bf \underline{Sub-case 2 (ii).}}\q $\lambda_n\geq \tilde{C}.$ 

Recall that by Assumption \ref{assumption_main} (ii), $\psi$ is differentiable in $[0,\t]$ and there exists $\k_2>0$ such that $|\psi'(s)|\leq \k_2$ on $[0,\t].$ Therefore, we have
\beqarray
\lambda_n^2 |\m_n| &=&\lambda_n^2\Big|\big[\frac{\psi(\t)}{\lambda_n^2}-\psi(0)(\frac{\t e^{-\lambda_n\t}}{\lambda_n}+\frac{e^{-\lambda_n\t}}{\lambda_n^2})\big]-\int_0^\t \psi_s(s)\big[\frac{(\t-s)e^{-\lambda_n(\t-s)}}{\lambda_n}+\frac{e^{-\lambda_n(\t-s)}}{\lambda_n^2}\big]\,ds\Big|\\
&\geq &|\psi(\t)|-\lambda_n^2|\psi(0)|\big[\frac{\t e^{-\lambda_n\t}}{\lambda_n}+\frac{e^{-\lambda_n\t}}{\lambda_n^2}\big]-\lambda_n^2\int_0^\t |\psi_s(s)|\big[\frac{(\t-s)e^{-\lambda_n(\t-s)}}{\lambda_n}+\frac{e^{-\lambda_n(\t-s)}}{\lambda_n^2}\big]\,ds\\
&\geq &|\psi(\t)|-\lambda_n^2|\psi(0)|\big[\frac{\t e^{-\lambda_n\t}}{\lambda_n}+\frac{e^{-\lambda_n\t}}{\lambda_n^2}\big]-\k_2\lambda_n^2\int_0^\t \big[\frac{(\t-s)e^{-\lambda_n(\t-s)}}{\lambda_n}+\frac{e^{-\lambda_n(\t-s)}}{\lambda_n^2}\big]\,ds\\
&=& |\psi(\t)|-\lambda_n^2|\psi(0)|\big[\frac{\lambda_n\t e^{-\lambda_n\t}+e^{-\lambda_n\t}}{\lambda_n^2}\big]-\lambda_n\k_2\big[ \frac{1}{\lambda_n^2}-\frac{\lambda_n\t e^{-\lambda_n\t}+e^{-\lambda_n\t}}{\lambda_n^2}\big]-\k_2\big[\frac{1}{\lambda_n}-\frac{e^{-\lambda_n\t}}{\lambda_n}\big]\\
&\geq & |\psi(\t)|-\lambda_n^2|\psi(0)|\big[\frac{\lambda_n\t e^{-\lambda_n\t}+e^{-\lambda_n\t}}{\lambda_n^2}\big]-\frac{2\k_2}{\lambda_n}\\
&\geq &|\psi(\t)|-\frac{|\psi(0)|}{\t\lambda_n}-\frac{2\k_2}{\lambda_n}\\
&\geq & |\psi(\t)|-\frac{|\psi(0)|}{\t\tilde{C}}-\frac{2\k_2}{\tilde{C}}\\
&=& |\psi(\t)|-\frac{1}{\tilde{C}}\big[\frac{|\psi(0)|+2\t{ \k_2}}{\t}\big].
\eeqarray
From this we obtain the required inequality with $C= |\psi(\t)|-\frac{1}{\tilde{C}}\big[\frac{|\psi(0)|+2\t{ \k_2}}{\t}\big].$
\epf

%\bcor
%Under Assumption \ref{assumption_main} we have $\int_0^\t e^{\lambda_n(s-\t)}(\t-s)\psi(s)\,ds\neq 0$ for all $n\in \N.$
%\ecor

Let $\varrho>0$ be fixed and $p> 0$. We now obtain a conditional stability estimate for the source set 
\beq\label{source_set}
S_{\varrho,\,p}=\{g\in L^2(\O): \sum_{n=1}^\infty \lambda_n^{2p}|\<g,\f_n\>|^2\leq \varrho^2\}=\{g\in \mathbb{H}_p: \|g\|_{\mathbb{H}_p}\leq \varrho\}.
\eeq

\bt\label{cond_stability_est}{\rm \bf (Conditional stability estimate)}
Let the exact source $f\in S_{\varrho,\,p}$, $h\in L^2(\O)$ and $C$ be the constant as in Theorem \ref{non_zero_lower_bound}. Then, we have $$\|f\|_{L^2(\O)}\leq C^{-\frac{p}{p+2}}\varrho^\frac{2}{p+2}\|h\|_{L^2(\O)}^\frac{p}{p+2}.$$
\et
\bpf
From the Fourier expansion of $f$, we have \beqarray
\|f\|^2_{L^2(\O)}&=&\sum_{n=1}^\infty |\<f,\f_n\>|^2=\sum_{n=1}^\infty \lambda_n^{\frac{4p}{p+2}}|\<f,\f_n\>|^\frac{4}{p+2}\lambda_n^{-\frac{4p}{p+2}}|\<f,\f_n\>|^\frac{2p}{p+2}\\
&\leq & \left(\sum_{n=1}^\infty \lambda_n^{2p}|\<f,\f_n\>|^2\right)^\frac{2}{p+2}\left(\sum_{n=1}^\infty \lambda_n^{-4}|\<f,\f_n\>|^2\right)^\frac{p}{p+2}\\
&\leq & \varrho^\frac{4}{p+2}\left(\sum_{n=1}^\infty \lambda_n^{-4}|\<f,\f_n\>|^2\right)^\frac{p}{p+2}.
\eeqarray
Hence, using (\ref{rep-h}) and Theorem \ref{non_zero_lower_bound}, we have 
$$\|f\|^2_{L^2(\O)} \leq \varrho^\frac{4}{p+2}\left(\sum_{n=1}^\infty \lambda_n^{-4} \Big|\frac{\<h,\f_n\>}{\m_n} \Big|^2\right)^\frac{p}{p+2}\\
\leq  \varrho^\frac{4}{p+2}C^{-\frac{2p}{p+2}}\|h\|_{L^2(\O)}^\frac{2p}{p+2}.$$ 
Hence, $\|f\|_{L^2(\O)}\leq C^{-\frac{p}{p+2}}\varrho^\frac{2}{p+2}\|h\|_{L^2(\O)}^\frac{p}{p+2}.$
\epf

\section{the quasi-reversibility method and error estimates}\label{sec-reg}

We have already observed that the  inverse source problem under consideration is ill-posed (also see, \cite[sec 3.1]{quoc_long_oregan_ngoc_tuan_2022}). We also illustrate this fact using an example. 
\beg{\bf (Ill-posed)}
Let $\tilde{h}=0$ and $\tilde{h}_k=\frac{1}{\lambda_k}\f_k$ for all $k\in \N.$ Let $\tilde{f}$ and $\tilde{f}_k$ denote the source functions corresponding to the final values $\tilde{h}$ and $\tilde{h}_k,$ respectively. Then by \eqref{rep_f_intermsof_h}, we have $\tilde{f}=0$ and $\tilde{f}_k=\frac{1}{\lambda_k\mu_k}\,\f_k.$ Thus, we have $\|\tilde{h}-\tilde{h}_k\|_{L^2(\O)}\to 0$ as $k\to \infty.$ But, from \eqref{bound_for_int_e_psi}, we have $\|\tilde{f}_k-\tilde{f}\|_{L^2(\O)}\geq \frac{1}{\|\psi\|_\infty}\,\lambda_k$ and hence $\|\tilde{f}_k-\tilde{f}\|_{L^2(\O)}\to \infty$ as $k\to \infty.$
\eeg
Thus, some  some regularization method has to be employed for obtaining stable approximate solutions. As mentioned earlier, we will be using the quasi-reversibility method (QRM) as follows. 

Let $b\geq 2$ and $\a>0.$ We consider the following system
\beq\label{quasi_rev_biparabolic_source}
\begin{cases}
v_{tt}-2\Delta v_t +\Delta^2v=(I+\a(-\Delta)^b)\psi(t)\tilde{f}(x),&\q\text{in}\q\O\times (0,\t),\\
v=0=\Delta v,&\q\text{on}\q\partial\O\times (0,\t),\\
v(\cdot,0)=0=v_t(\cdot,0),&\q\text{in}\q\O,\\
v(\cdot,\t)=\tilde{h},&\q\text{in}\q\O.
\end{cases}
\eeq
Suppose $\psi\in L^2(0,\t),\,\,\tilde{h}\in L^2(\O)$ and $\tilde{f}\in \mathbb{H}_b$, the domain of $(-\Delta)^b$. Suppose $v$ is a solution of \eqref{quasi_rev_biparabolic_source}.  Then following \cite{quoc_long_oregan_ngoc_tuan_2022, tuan_2022} we observe the following:

Taking the $L^2$ inner product on both sides of the governing equation of \eqref{quasi_rev_biparabolic_source} with $\f_n$, for all $n\in \N,$ we obtain
$$\psi(t)\<(I+\a(-\Delta)^b)\tilde{f},\f_n\>=\frac{d^2}{dt^2}\< v(t),\f_n\>-2\frac{d}{dt}\< \Delta v(t),\f_n\>+\<\Delta^2 v(t),\f_n\>.$$
Now for $\tilde{f}\in \mathbb{H}_b,$ we have $\ds(I+(-\Delta)^b)\tilde{f}=\sum_{n=1}^\infty (1+\a\lambda_n^b)\<\tilde{f},\f_n\>\f_n.$ Therefore, 
$$\psi(t)\<\tilde{f},\f_n\>=\frac{1}{1+\a\lambda_n^b}\left(\frac{d^2}{dt^2}\<v(t),\f_n\>-2\frac{d}{dt}\<\Delta v(t),\f_n\>+\<\Delta^2 v(t),\f_n\>\right)$$
Now, doing repeated integration by parts and making use of the boundary condition in \eqref{quasi_rev_biparabolic_source}, we obtain
$$
\<\Delta v(t),\f_n\>=\<v(t),\Delta \f_n\>=-\lambda_n\<v(t),\f_n\>
$$
and
$$\<\Delta^2 v(t),\f_n\>=\<v(t),\Delta^2 \f_n\>=\lambda_n^2\<v(t),\f_n\>.$$
Therefore, we have
\beq\label{reduced_ODE}
\psi(t)\<\tilde{f},\f_n\>=\frac{1}{1+\a\lambda_n^b}\left(\frac{d^2}{dt^2}\<v(t),\f_n\>+2\lambda_n\frac{d}{dt}\<v(t),\f_n\>+\lambda_n^2\<v(t),\f_n\>\right).
\eeq
Now it is easy to check that 
$$\<v(t),\f_n\> = (1+\a\lambda_n^b) \mu_n(t)\<\tilde{f},\f_n\>,$$
where $\m_n(t)$ is as in (\ref{ev-t}). Hence,  the solution $v$ is given by
\beq\label{sol_qrm_smooth} 
v(t):=\sum_{n=1}^\infty  (1+\a\lambda_n^b) \mu_n(t)\<\tilde{f},\f_n\>\f_n,\eeq

Now, analogous to the observation \eqref{bound_for_int_e_psi}, we have 
$$ |\m_n(t)| = \Big|\int_0^t e^{-\lambda_n(t-s)}(t-s)\psi(s)\,ds\Big|\leq \frac{\|\psi\|_\infty}{\lambda_n^2}.$$ Thus, in order to ensure $v(t)\in L^2(\O)$ we need $\tilde{f}\in \mathbb{H}_{b-2}$ only, and not necessarily in $\mathbb{H}_b,$ the domain of $(-\Delta)^b$. %In view of this, we define the {\it mild solution} for \eqref{quasi_rev_biparabolic_source}.

In view of the above discussion, we introduce the following definition. 

\bedef\label{mild_sol_qrm}
For $\psi\in L^2(0,\t)$ and $\tilde{f}\in \mathbb{H}_{b-2},$ a function $v\in C([0,\t];L^2(\O))$ with the representation given by \eqref{sol_qrm_smooth} is said to be a mild solution of \eqref{quasi_rev_biparabolic_source}.
\endef

%%%%%%%  Remark 4.3 has been removed %%%%

From the above discussion, we observe that for $\tilde{h}\in L^2(\O)$ and $\tilde{f}\in \mathbb{H}_{b-2},\,b\geq 2,$ \eqref{quasi_rev_biparabolic_source} has a unique mild solution, and we call it as $v.$
Now, taking into account the final value condition $v(\cdot,\t)=\tilde{h}$, we obtain
$$ \<\tilde{f},\f_n\>=\frac{ \<\tilde{h},\f_n\>}{(1+\a\lambda_n^b) \m_n}.$$
Thus, the problem of identification of $\tilde{f}$ from the final value $\tilde{h}$ associated with \eqref{quasi_rev_biparabolic_source} is same as the problem of solving the operator equation
\beq\label{qrm_opeq}\mathbb{T}_\a \tilde{f}=\tilde{h},\eeq
where $\mathbb{T}_\a:\text{\bf dom}(\mathbb{T}_\a)\subset L^2(\O)\to L^2(\O)$ is defined by
$$\mathbb{T}_\a\phi:=\sum_{n=1}^\infty (1+\a\lambda_n^b)\mu_n\<\phi,\f_n\>\f_n.$$
Since $|\mu_n|\sim \frac{1}{\lambda_n^2}$, it is easy to check that $\text{\bf dom}(\mathbb{T}_\a)=\mathbb{H}_{b-2}.$ Also, it is easy to check that $\mathbb{T}_\a$ is a densely defined closed self-adjoint operator. Now, for $b\geq 2,\,\,(1+\a\lambda_n^b)^b\geq (1+\a\lambda_n^b)^2\geq\a^2\lambda_n^{2b}$ and hence $(1+\a\lambda_n^b)\geq \a^\frac{2}{b}\lambda_n^2$. In view of this and from Theorem \ref{non_zero_lower_bound}, we have
$$\|\mathbb{T}_\a\phi\|_{L^2(\O)}\geq C\a^\frac{2}{b}\|\phi\|_{L^2(\O)}\q\text{for all}\q\phi\in \text{\bf dom}(\mathbb{T}_\a).$$
Thus, $\mathbb{T}_\a^{-1}:\text{\bf ran}(\mathbb{T}_\a)\to L^2(\O)$ is a closed operator which is also bounded, and hence $\text{\bf ran}(\mathbb{T}_\a)$ is closed (cf. \cite{nair_opeq}). Now from the fact $\text{\bf ran}(\mathbb{T}_\a)=N(\mathbb{T}_\a)^\perp=\{0\}^\perp$, we have $\text{\bf ran}(\mathbb{T}_\a)=L^2(\O).$ Therefore, it follows that $\mathbb{T}_\a$ has a bounded inverse with $$\|\mathbb{T}_\a^{-1}\phi\|_{L^2(\O)}\leq \frac{1}{C\a^\frac{2}{b}}\|\phi\|_{L^2(\O)}\q\text{for all}\q\phi\in L^2(\O).$$

Thus, from the above discussion, we observe that solving equation \eqref{qrm_opeq} is an well-posed problem. Moreover, for $\a>0$, let $f_\a:=\mathbb{T}_\a^{-1}h$ and $f_\a^\d:=\mathbb{T}_\a^{-1}h^\d.$ Then, we have
$$\|f_\a-f_\a^\d\|_{L^2(\O)}\leq \frac{1}{C\a^\frac{2}{b}}\|h-h^\d\|_{L^2(\O)}.$$

%\bt{\bf(Well-posedness of solution by QRM)}\label{qrm_well_posed}
%Let $\tilde{h}\in L^2(\O)$, $b\geq 2$, $\a>0,$ and $C$ be as in Theorem \ref{non_zero_lower_bound}. Then \eqref{qrm_opeq} has a unique solution $\tilde{f}:=\sum_{n=1}^\infty \frac{\<\tilde{h},\f_n\>}{(1+\a\lambda_n^b)\mu_n}\f_n.$ Moreover, $\|\tilde{f}\|_{L^2(\O)}\leq \frac{1}{C\a^\frac{2}{b}}\|h\|_{L^2(\O)},$ that is, $\tilde{f}$ depends continuously on $\tilde{h}.$
%\et
%\bpf
%The existence of the solution is clear from the definition of $\mathbb{T}_\a.$ We only show the continuous dependence. From the representation of $\tilde{f},$ we have
%$$\|\tilde{f}\|^2_{L^2(\O)}=\sum_{n=1}^\infty \Big|\frac{\<\tilde{h},\f_n\>}{(1+\a\lambda_n^b) \m_n }\Big|^2 \leq  \frac{1}{\a^\frac{4}{b}}\sum_{n=1}^\infty \Big|\frac{\<\tilde{h},\f_n\>}{\lambda_n^2 \m_n }\Big|^2
%\underbrace{\leq }_{{\rm Theorem} \ref{non_zero_lower_bound}}  \frac{1}{C^2\,\a^\frac{4}{b}}\|h\|^2_{L^2(\O)}$$
%\epf

%This shows that the inverse problem of source identification associated with \eqref{quasi_rev_biparabolic_source} has a unique solution for the final value. 

Recall that, $f_\a=\mathbb{T}_\a^{-1}h=\sum_{n=1}^\infty\frac{ \<h,\f_n\>}{(1+\a\lambda_n^b) \m_n }\f_n$ and $f_\a^\d=\mathbb{T}_\a^{-1}h^\d=\sum_{n=1}^\infty\frac{ \<h^\d,\f_n\>}{(1+\a\lambda_n^b) \m_n}\f_n.$  Therefore, we have
\beq\label{f_alpha_and_f_alpha_delta_rep_interms_h_and_hdelta}
\<f_\a,\f_n\>=\frac{ \<h,\f_n\>}{(1+\a\lambda_n^b) \m_n }\q\text{and}\q \<f_\a^\d,\f_n\>=\frac{ \<h^\d,\f_n\>}{(1+\a\lambda_n^b) \m_n}.
\eeq
Before proceeding further, let us mention that Young's inequality will be used in many places in the upcoming analysis without mentioning explicitly. The Young's inequality states that: for $a,b\geq 0$ and $\theta, \theta'\geq 1$ satisfying $\frac{1}{\theta}+\frac{1}{\theta'}=1$, we have $ab\leq \frac{a^\theta}{\theta}+\frac{b^{\theta'}}{\theta'}.$
\subsection{Apriori parameter choice}
We are now in a position to state and prove one of the main results of this paper, the order of the error estimate for an apriori parameter choice strategy. First we derive some general error estimates in terms of $(\a,\d)$.
\bt\label{apriori_preparatory}
Let $\a,\d>0, b\geq 2$ and $h^\d\in L^2(\O)$ be as in \eqref{noise_model}. Let $f\in S_{\varrho,\,p}$ for some $\varrho, p>0$ and $C$ be as in Theorem \ref{non_zero_lower_bound}. Then we have the following estimates:
\ben
\i[(i)] $\|f_\a-f_\a^\d\|_{L^2(\O)}\leq \frac{1}{C}\frac{\d}{\a^\frac{2}{b}}.$
\i[(ii)] $\|f-f_\a\|_{L^2(\O)}\leq 
\begin{cases}
C_{\rm apri1}\,\varrho\,\a^\frac{p}{b},&\q\text{if}\q 0<p<b,\\
C_{\rm apri2}\,\varrho\,\a,&\q\text{if}\q p\geq b,
\end{cases}$
\newline where $C_{\rm apri1}=\frac{p}{b}\,\big(\frac{b-p}{p}\big)^\frac{b-p}{b}$ and $C_{\rm apri2}=\lambda_1^{(b-p)}.$
\een
\et
\bpf
{\bf  (i).} First, we observe that for $b\geq 2,$
$$(1+\a\lambda_n^b)^b\geq (1+\a\lambda_n^b)^2\geq\a^2\lambda_n^{2b}$$
and hence $1+\a\lambda_n^b\geq \a^\frac{2}{b}\lambda_n^2.$
%\beqarray \cred 
%1+\a\lambda_n^b\geq \frac{b-2}{b} + \frac{2}{b}\,\a\lambda_n^b=\frac{b-2}{b}\,[1^\frac{b-2}{b}]^\frac{b}{b-2}+\frac{2}{b}[(\a\lambda_n^b)^\frac{2}{b}]^\frac{b}{2}\geq (\a\lambda_n^b)^\frac{2}{b}=\a^\frac{2}{b}\lambda_n^2.
%\eeqarray
Therefore, from \eqref{f_alpha_and_f_alpha_delta_rep_interms_h_and_hdelta} it follows that
 
\beqarray 
\|f_\a^\d-f_\a\|^2_{L^2(\O)}&=&\sum_{n=1}^\infty \Big|\frac{\<h^\d-h,\f_n\>}{(1+\a\lambda_n^b) \m_n }\Big|^2 \leq \frac{1}{\a^\frac{4}{b}}\sum_{n=1}^\infty \Big|\frac{\<h^\d-h,\f_n\>}{\lambda_n^2 \m_n }\Big|^2\\
&\underbrace{\leq }_{\rm Theorem \ref{non_zero_lower_bound}}&\frac{1}{\a^\frac{4}{b}C^2}\sum_{n=1}^\infty |\<h^\d-h,\f_n\>|^2 \leq \frac{1}{C^2}\frac{\d^2}{\a^\frac{4}{b}}.
\eeqarray
Hence, we obtain
\beq\label{falphadelta_falpha_est}
\|f_\a^\d-f_\a\|_{L^2(\O)}\leq \frac{1}{C}\frac{\d}{\a^\frac{2}{b}}.
\eeq

\noi
{\bf (ii).} From \eqref{f_alpha_and_f_alpha_delta_rep_interms_h_and_hdelta} we have 
\beqarray
\|f-f_\a\|^2_{L^2(\O)}&=&\sum_{n=1}^\infty \Big| \frac{\<h,\f_n\>}{\m_n }-\frac{\<h,\f_n\>}{(1+\a\lambda_n^b) \m_n }\Big|^2 = \sum_{n=1}^\infty \Big|\frac{\a\lambda_n^b\<h,\f_n\>}{(1+\a\lambda_n^b) \m_n }\Big|^2\\
&=& \sum_{n=1}^\infty \frac{\a^2\lambda_n^{2b}}{(1+\a\lambda_n^b)^2}|\<f,\f_n\>|^2= \sum_{n=1}^\infty \frac{\a^2\lambda_n^{2b-2p}}{(1+\a\lambda_n^b)^2}\lambda_n^{2p}|\<f,\f_n\>|^2.
\eeqarray
Now, we consider two cases, namely, $0<p<b$ and $p\geq b$.
\vsq 
\noi {\bf \underline{Case 1.}} $0<p<b$\,:  
In this case, by Lemma \ref{lemma-basic}(ii),  we have $$\frac{\a^2\lambda_n^{2b-2p}}{(1+\a\lambda_n^b)^2}\leq \left(\frac{p}{b}\right)^2\left(\frac{b-p}{p}\right)^\frac{2(b-p)}{b}\a^\frac{2p}{b}$$ and hence
\beq\label{apriori_est_f_and_falpha_case1}
\|f-f_\a\|_{L^2(\O)}\leq C_{\rm apri1}\,\varrho\, \a^\frac{p}{b},
\eeq
where $C_{\rm apri1}=\frac{p}{b}\,\left(\frac{b-p}{p}\right)^\frac{b-p}{b}.$

\noi {\bf \underline{Case 2.}}  $p\geq b$\,: 
In this case, $$\frac{\a^2\lambda_n^{2b-2p}}{(1+\a\lambda_n^b)^2}\leq \a^2\lambda_1^{2b-2p},$$
and hence 
\beq\label{apriori_est_f_and_falpha_case2}
\|f-f_\a\|_{L^2(\O)}\leq C_{\rm apri2}\,\varrho\,\a,
\eeq
where $C_{\rm apri2}=\lambda_1^{(b-p)}.$
\epf
Now we are in a position to state the main result of this subsection.
\bt\label{qrm_err_est_apriori}
Let $\a,\d>0, b\geq 2$ and $h^\d\in L^2(\O)$ be as in \eqref{noise_model}. Let $f\in S_{\varrho,\,p}$ for some $\varrho, p>0.$ Then the following holds.
\ben
\i[(i)] If $0<p<b$ then for $\a\sim \left(\frac{\d}{\varrho}\right)^\frac{b}{p+2}$ there exists a constant $C_1:=C_1(C,p,b)>0$ such that $$\|f-f_\a^\d\|_{L^2(\O)}\leq C_1\,\varrho^\frac{2}{p+2}\,\d^\frac{p}{p+2}.$$
\i[(ii)] If $p\geq b$ then for $\a\sim \left(\frac{\d}{\varrho}\right)^\frac{b}{b+2}$ there exists a constant $C_2:=C_2(C,\lambda_1,p,b)>0$ such that
$$\|f-f_\a^\d\|_{L^2(\O)}\leq C_2\,\varrho^\frac{2}{b+2}\,\d^\frac{b}{b+2}.$$
\een
\et
\bpf
By Theorem \ref{apriori_preparatory}, we have
$$\|f-f_\a^\d\|_{L^2(\O)}\leq 
\begin{cases}
C_{\rm apri1}\,\varrho\,\a^\frac{p}{b}+\frac{1}{C}\,\frac{\d}{\a^\frac{2}{b}},&\q\text{if}\q 0<p<b,\\
C_{\rm apri2}\,\varrho\,\a+\frac{1}{C}\frac{\d}{\a^\frac{2}{b}},&\q\text{if}\q p\geq b.
\end{cases}$$
Hence the proof follows by appropriately choosing the constants $C_1$ and $C_2.$
\epf
\brem\label{better_rate_apriori}
If $b=2,$ then by Theorem \ref{qrm_err_est_apriori}, the best possible order of the error estimate in the apriori case is $\mathcal{O}(\d^\frac{1}{2})$ which is better than the order $\mathcal{O}(\d^\frac{1}{3})$ obtained in \cite[Theorem 4.1]{quoc_long_oregan_ngoc_tuan_2022} using Tikhonov regularization method and the order $\mathcal{O}(\d^\frac{1}{4})$ obtained in \cite[Theorem 3.2]{phuong_luc_long_2020} using modified quasi-boundary value method for the apriori case.

Also, it is easy to observe that if $b\geq 4$ then $\frac{b}{b+2}\geq \frac{2}{3}$ and hence by Theorem \ref{qrm_err_est_apriori} it follows that the best possible order can exceed the benchmark order of $\mathcal{O}(\d^\frac{2}{3})$ for Tikhonov regularization. 
\erem

\subsection{Aposteriori parameter choice}
We now proceed on to aposteriori case, where we will do the analysis using Morozov-type discrepancy principle, see for e.g., \cite{engl_hanke_neubauer, nair_opeq}. For $\a>0,$ let $\mathbb{B}_\a:L^2(\O)\to L^2(\O)$ be defined by $$\mathbb{B}_\a\tilde{h}:=\sum_{n=1}^\infty \frac{\<\tilde{h},\f_n\>}{(1+\a\lambda_n^b)}\f_n.$$  Then we have the following result.
\bl\label{apost_zeta}(cf. \cite{hao_duc_thang_thanh_2020})
For $\a>0,$ let $\zeta(\a):=\|\mathbb{B}_\a h^\d-h^\d\|_{L^2(\O)}.$ Then 
\ben
\i[(i)] $\zeta$ is a continuous function.
\i[(ii)] $\lim\limits_{\a\to 0^+}\zeta(\a)=0$ and $\lim\limits_{\a\to +\infty}\zeta(\a)=h^\d.$
\i[(iii)] $\zeta$ is strictly increasing.
\een
\el

We are now in a position to state another main result of this paper, namely the error estimates for an aposteriori parameter choice strategies.
\bt\label{qrm_err_est_apost}
Let the Assumption \ref{assumption_main} holds, $f\in S_{\varrho,p}$, $\d>0$ and $h^\d$ be as in \eqref{noise_model}. Let $\s\in (0,1)$ and $b\geq 2.$ Suppose that $\d^\s<\|h^\d\|_{L^2(\O)}.$ Choose $\xi>1$ such that $0<\xi \d^\s\leq \|h^\d\|_{L^2(\O)}.$ Then we have the following estimates.
\ben
\i[(i)] If $b\neq 2$ then for each $\d>0$ there exists a unique $\a_\d>0$ such that $\|\mathbb{B}_{\a_\d}h^\d-h^\d\|_{L^2(\O)}=\xi\d.$ Further, there exists constants  $C_4:=C_4(C,\|\psi\|_\infty,\xi,\lambda_1,p,b), C_6:=C_6(C,\|\psi\|_\infty,\xi,p), C_{\rm apost}:=C_{\rm apost}(C,\xi,p)>0$ such that 
$$
\|f-f_{\a_\d}^\d\|_{L^2(\O)}\leq \begin{cases}
C_6\,\varrho^\frac{2}{p+2}\,\d^\frac{p}{p+2},&\q\text{if}\q 0<p<b-2\\
C_{\rm apost}\,\varrho^\frac{2}{p+2}\,\d^\frac{p}{p+2}+C_4\,\varrho^\frac{2}{b}\,\d^\frac{b-2}{b},&\q\text{if}\q p\geq b-2.
\end{cases}$$
\i[(ii)] If $b= 2$ then for each $\d>0$ there exists a unique $\a_\d>0$ such that $\|\mathbb{B}_{\a_\d}h^\d-h^\d\|_{L^2(\O)}=\xi\d^\s.$ Further, there exists constant  $C_5:=C_5(C,\|\psi\|_\infty,\xi,\lambda_1,p)>0$ (and $C_{\rm apost}$ as in (i)) such that
$$\|f-f_{\a_\d}^\d\|_{L^2(\O)}\leq C_{\rm apost}\,\varrho^\frac{2}{p+2}\,\d^\frac{p\s}{p+2}+C_5\,\varrho\,\d^{1-\s}.$$
\een
\et
\bpf
Let $w$ be the mild solution (see Definition \eqref{mild_sol_qrm}) of 
$$
\begin{cases}
v_{tt}-2\Delta v_t +\Delta^2v=(I+\a(-\Delta)^b)\psi(t)f_\a(x),&\q\text{in}\q\O\times (0,\t),\\
v=0=\Delta v,&\q\text{on}\q\partial\O\times (0,\t),\\
v(\cdot,0)=0=v_t(\cdot,0),&\q\text{in}\q\O,\\
v(\cdot,\t)=h,&\q\text{in}\q\O.
\end{cases}
$$
Then it is easy to check that $z:=\mathbb{B}_\a w$ is a solution of
\beq
\begin{cases}
v_{tt}-2\Delta v_t +\Delta^2v=\psi(t)f_\a(x),&\q\text{in}\q\O\times (0,\t),\\
v=0=\Delta v,&\q\text{on}\q\partial\O\times (0,\t),\\
v(\cdot,0)=0=v_t(\cdot,0),&\q\text{in}\q\O,\\
v(\cdot,\t)=\mathbb{B}_\a h,&\q\text{in}\q\O.
\end{cases}
\eeq
Now,
$$ \|f-f_\a\|^2_{\mathbb{H}_p}=\sum_{n=1}^\infty \lambda_n^{2p}\big[\<f,\f_n\>-\frac{\<h,\f_n\>}{(1+\a_\d\lambda_n^b) \m_n }\big]^2=\sum_{n=1}^\infty \lambda_n^{2p}\big[\frac{\a_\d\lambda_n^b\<f,\f_n\>}{1+\a_\d\lambda_n^b}\big]^2\leq \|f\|^2_{\mathbb{H}_p}\leq \varrho^2$$
and hence
\beq\label{apost_stability_est}
\|f-f_\a\|_{\mathbb{H}_p}\leq \varrho.
\eeq

\noi{\bf \underline{Proof of (i).}}\q By Lemma \ref{apost_zeta} it follows that there  exists a unique $\a_\d>0$ such that $\|\mathbb{B}_{\a_\d}h^\d-h^\d\|_{L^2(\O)}=\xi \d.$ 
Also, we have
\beqarray
\|u(\t)-z(\t)\|_{L^2(\O)}&=&\|h-\mathbb{B}_{\a_\d} h\|_{L^2(\O)}\leq \|\mathbb{B}_{\a_\d}h-\mathbb{B}_{\a_\d}h^\d\|_{L^2(\O)}+\|\mathbb{B}_{\a_\d}h^\d-h^\d\|_{L^2(\O)}+\|h^\d-h\|_{L^2(\O)}\\
&\leq & \|\sum_{n=1}^\infty \frac{\<h-h^\d,\f_n\>}{(1+\a\lambda_n^b)}\f_n\|_{L^2(\O)}+\xi \d+\d\leq (\xi+2)\d.
\eeqarray
Note that $u-z$ is a solution of 
$$
\begin{cases}
v_{tt}-2\Delta v_t +\Delta^2v=\psi(t)\,(f-f_{\a_\d})(x),&\q\text{in}\q\O\times (0,\t),\\
v=0=\Delta v,&\q\text{on}\q\partial\O\times (0,\t),\\
v(\cdot,0)=0=v_t(\cdot,0),&\q\text{in}\q\O,\\
v(\cdot,\t)=h-\mathbb{B}_{\a_\d} h,&\q\text{in}\q\O.
\end{cases}
$$ 
and \eqref{apost_stability_est} shows that $f-f_{\a_\d}\in S_{\varrho, p}.$ Therefore, following the proof of Theorem \ref{cond_stability_est}, it is easy to check that
$\|f-f_{\a_\d}\|_{L^2(\O)}\leq C^{-\frac{p}{p+2}}\varrho^\frac{2}{p+2}[(\xi+2)\d]^\frac{p}{p+2}$ and hence
\beq\label{est_apost_f_and_falpha_apost1}
\|f-f_{\a_\d}\|_{L^2(\O)}\leq C_{\rm apost}\,\varrho^\frac{2}{p+2}\,\d^\frac{p}{p+2},
\eeq
where $C_{\rm apost}=\left(\frac{\xi+2}{C}\right)^\frac{p}{p+2}.$ We now try to obtain a bound for $\a_\d$ in terms of $\d$ using the discrepancy principle. Now,
\beqarray
\xi\d&=&\|\mathbb{B}_{\a_\d}h^\d-h^\d\|_{L^2(\O)}=\|\sum_{n=1}^\infty \frac{\a_\d\lambda_n^b}{1+\a_\d\lambda_n^b}\<h^\d,\f_n\>\f_n\|_{L^2(\O)}\\
&\leq & \|\sum_{n=1}^\infty \frac{\a_\d\lambda_n^b}{1+\a_\d\lambda_n^b}\<h,\f_n\>\f_n\|_{L^2(\O)}+\|\sum_{n=1}^\infty \frac{\a_\d\lambda_n^b}{1+\a_\d\lambda_n^b}\<h^\d-h,\f_n\>\f_n\|_{L^2(\O)}\\
&\leq & \|\sum_{n=1}^\infty \frac{\a_\d\lambda_n^b}{1+\a_\d\lambda_n^b} \m_n \,\<f,\f_n\>\f_n\|_{L^2(\O)}+\d.\eeqarray 
Therefore, 
$$ (\xi-1)\d \leq \|\sum_{n=1}^\infty \frac{\a_\d\lambda_n^b}{1+\a_\d\lambda_n^b} \m_n \,\<f,\f_n\>\f_n\|_{L^2(\O)}.$$
Thus,
\beq\label{apost_aux_est}
(\xi-1)^2\d^2\leq \sum_{n=1}^\infty \left(\frac{\a_\d\lambda_n^{b-2}}{1+\a_\d\lambda_n^b}\right)^2\big(\<f,\f_n\>\|\psi\|_{\infty}\big)^2
\eeq
{\bf \underline{Case 1.}}\q $b>2$ and $0<p<b-2.$

In this case,
\beqarray
\a_\d\lambda_n^b+1&\geq & \frac{b-2-p}{b}\a_\d\lambda_n^b+\frac{p+2}{b}\\
&=& \frac{b-2-p}{b} \big[(\a_\d\lambda_n^b)^\frac{b-2-p}{b}\big]^\frac{b}{b-2-p}+\frac{p+2}{b}[1^\frac{p+2}{b}]^\frac{b}{p+2}\\
&\geq & \left(\a_\d\lambda_n^b\right)^\frac{b-2-p}{b},
\eeqarray
and hence $\frac{\a_\d\lambda_n^{b-2}}{1+\a_\d\lambda_n^b}\leq \a_\d^\frac{p+2}{b}\lambda_n^p.$
Therefore,
\beqarray
&&[(\xi-1)\d]^2\leq \|\psi\|^2_{\infty}\,\a_\d^{\frac{2(p+2)}{b}}\sum_{n=1}^\infty \lambda_n^{2p}|\<f,\f_n\>|^2\leq \|\psi\|^2_\infty\,\a_\d^{\frac{2(p+2)}{b}}\,\varrho^2\\
&\implies & \frac{1}{\a_\d}\leq \left(\frac{\|\psi\|_\infty}{\xi-1}\right)^\frac{b}{p+2} \left(\frac{\varrho}{\d}\right)^\frac{b}{p+2}.
\eeqarray
Thus, from \eqref{falphadelta_falpha_est}, we have
$$
\|f_{\a_\d}^\d-f_{\a_\d}\|_{L^2(\O)}\leq \frac{1}{C}\frac{\d}{{\a_\d}^\frac{2}{b}}\leq \frac{1}{C}\,\d\,\left(\frac{\|\psi\|_\infty}{\xi-1}\right)^\frac{2}{p+2}\left(\frac{\varrho}{\d}\right)^\frac{2}{p+2}=\frac{\|\psi\|^\frac{2}{p+2}_\infty}{C(\xi-1)^\frac{2}{p+2}}\,\varrho^\frac{2}{p+2}\,\d^\frac{p}{p+2},
$$
that is,
\beq\label{apost_est_falpha_delta_falpha_case11}
\|f_{\a_\d}^\d-f_{\a_\d}\|_{L^2(\O)}\leq C_3\,\varrho^\frac{2}{p+2}\,\d^\frac{p}{p+2},
\eeq
where $C_3=\frac{\|\psi\|^\frac{2}{p+2}_\infty}{C(\xi-1)^\frac{2}{p+2}}.$

\noi {\bf \underline{Case 2.}}\q $b>2$ and $p\geq b-2.$

In this case, from \eqref{apost_aux_est} we have,
\beqarray
&&[(\xi-1)\d]^2\leq \|\psi\|^2_\infty\sum_{n=1}^\infty  \Big[\frac{\a_\d\lambda_n^{b-2-p}}{1+\a_\d\lambda_n^b}\lambda_n^p\<f,\f_n\>\Big]^2\leq \frac{\|\psi\|^2_\infty}{\lambda_1^{2(p-b+2)}}\,\a_\d^2\,\varrho^2,\\
&\implies & \frac{1}{\a_\d}\leq \frac{\|\psi\|_\infty}{(\xi-1)\lambda_1^{p-b+2}}\,\frac{\varrho}{\d}.
\eeqarray
Therefore, from \eqref{falphadelta_falpha_est}, we have
$$
\|f_{\a_\d}^\d-f_{\a_\d}\|_{L^2(\O)}\leq \frac{1}{C}\frac{\d}{{\a_\d}^\frac{2}{b}}\leq \frac{1}{C}\,\d\,\Big[\frac{\|\psi\|_\infty}{(\xi-1)\lambda_1^{p-b+2}}\Big]^\frac{2}{b}\,\left(\frac{\varrho}{\d}\right)^\frac{2}{b}=\frac{1}{C} \Big[\frac{\|\psi\|_\infty}{(\xi-1)\lambda_1^{p-b+2}}\Big]^\frac{2}{b}\,\varrho^\frac{2}{b}\,\d^\frac{b-2}{b}
$$
that is,
\beq\label{apost_est_falpha_delta_falpha_case12}
\|f_{\a_\d}^\d-f_{\a_\d}\|_{L^2(\O)}\leq C_4\,\varrho^\frac{2}{b}\,\d^\frac{b-2}{b},
\eeq
where $C_4=\frac{1}{C} \Big[\frac{\|\psi\|_\infty}{(\xi-1)\lambda_1^{p-b+2}}\Big]^\frac{2}{b}.$ 

Thus, from \eqref{est_apost_f_and_falpha_apost1}, \eqref{apost_est_falpha_delta_falpha_case11} and \eqref{apost_est_falpha_delta_falpha_case12} we have
$$
\|f-f_{\a_\d}^\d\|_{L^2(\O)}\leq \begin{cases}
C_{\rm apost}\,\varrho^\frac{2}{p+2}\,\d^\frac{p}{p+2}+C_3\,\varrho^\frac{2}{p+2}\,\d^\frac{p}{p+2},&\q\text{if}\q 0<p<b-2\\
C_{\rm apost}\,\varrho^\frac{2}{p+2}\,\d^\frac{p}{p+2}+C_4\,\varrho^\frac{2}{b}\,\d^\frac{b-2}{b},&\q\text{if}\q p\geq b-2.
\end{cases}$$
Now the proof follows by choosing the constant $C_6$ appropriately.

\noi {\bf \underline{Proof of (ii).}}\q Again by Lemma \ref{apost_zeta}, it follows that there exists a unique $\a_\d>0$ such that $\|\mathbb{B}_{\a_\d}(h^\d)-h^\d\|_{L^2(\O)}=\xi \d^\s$ for $\s\in (0,1).$
Now, we have
\beqarray
\|u(\t)-z(\t)\|_{L^2(\O)}&=&\|h-\mathbb{B}_{\a_\d} h\|_{L^2(\O)}\leq \|\mathbb{B}_{\a_\d}h-\mathbb{B}_{\a_\d}h^\d\|_{L^2(\O)}+\|\mathbb{B}_{\a_\d}h^\d-h^\d\|_{L^2(\O)}+\|h^\d-h\|_{L^2(\O)}\\
&\leq & \d+\xi\d^\s+\d\leq (\xi+2)\d^\s.
\eeqarray
Therefore from \eqref{apost_stability_est} and Theorem \ref{cond_stability_est}, we have
$\|f-f_{\a_\d}\|_{L^2(\O)}\leq C^{-\frac{p}{p+2}}\varrho^\frac{2}{p+2}[(\xi+2)\d^\s]^\frac{p}{p+2}$ and hence
\beq\label{est_apost_f_and_falpha_apost2}
\|f-f_{\a_\d}\|_{L^2(\O)}\leq C_{\rm apost}\,\varrho^\frac{2}{p+2}\,\d^\frac{p\s}{p+2},
\eeq
where $C_{\rm apost}=\Big(\frac{\xi+2}{C}\Big)^\frac{p}{p+2}.$

Now, similar to {\bf Case 1}, we obtain a bound for $\a_\d$ in terms of $\d$ using the discrepancy principle. Now,
\beqarray
\xi\d^\s&=&\|\mathbb{B}_{\a_\d}h^\d-h^\d\|_{L^2(\O)}=\|\sum_{n=1}^\infty \frac{\a_\d\lambda_n^2}{1+\a_\d\lambda_n^2}\<h^\d,\f_n\>\f_n\|_{L^2(\O)}\\
%&\leq & \|\sum_{n=1}^\infty \frac{\a_\d\lambda_n^2}{1+\a_\d\lambda_n^2}\<h,\f_n\>\f_n\|_{L^2(\O)}+\|\sum_{n=1}^\infty \frac{\a_\d\lambda_n^2}{1+\a_\d\lambda_n^2}\<h^\d-h,\f_n\>\f_n\|_{L^2(\O)}\\
&\leq & \|\sum_{n=1}^\infty \frac{\a_\d\lambda_n^2}{1+\a_\d\lambda_n^2} \m_n \,\<f,\f_n\>\f_n\|_{L^2(\O)}+\d 
\eeqarray 
so that  
$$ (\xi-1)\d^\s \leq \|\sum_{n=1}^\infty \frac{\a_\d\lambda_n^2}{1+\a_\d\lambda_n^2} \m_n \,\<f,\f_n\>\f_n\|_{L^2(\O)}.$$
Thus, by \eqref{bound_for_int_e_psi}, we have
\beqarray
&&(\xi-1)^2\d^{2\s}\leq \sum_{n=1}^\infty \Big[\frac{\a_\d}{1+\a_\d\lambda_n^2}\<f,\f_n\>\|\psi\|_\infty\Big]^2\leq \frac{\|\psi\|^2_\infty}{\lambda_1^{2p}}\,\varrho^2\,\a_\d^2\\
&\implies & \frac{1}{\a_\d}\leq \frac{\|\psi\|_\infty}{(\xi-1)\lambda_1^p }\,\frac{\varrho}{\d^\s}.
\eeqarray
Hence, from \eqref{falphadelta_falpha_est}, we have
$$\|f_{\a_\d}^\d-f_{\a_\d}\|_{L^2(\O)}\leq \frac{1}{C}\frac{\d}{\a_\d}\leq \frac{1}{C}\,\d\,\frac{\|\psi\|_\infty}{(\xi-1)\lambda_1^p } \,\frac{\varrho}{\d^\s}=\frac{1}{C}\frac{\|\psi\|_\infty}{(\xi-1)\lambda_1^p }\,\varrho\,\d^{1-\s},$$
that is,
\beq\label{apost_est_falpha_delta_falpha_case2}
\|f_{\a_\d}^\d-f_{\a_\d}\|_{L^2(\O)}\leq C_5\,\varrho\,\d^{1-\s},
\eeq
where $C_5=\frac{1}{C}\frac{\|\psi\|_\infty}{(\xi-1)\lambda_1^p }.$

Thus, from \eqref{est_apost_f_and_falpha_apost2} and \eqref{apost_est_falpha_delta_falpha_case2}, we have
$$\|f-f_{\a_\d}^\d\|_{L^2(\O)}\leq C_{\rm apost}\,\varrho^\frac{2}{p+2}\,\d^\frac{p\s}{p+2}+C_5\,\varrho\,\d^{1-\s}.$$
This completes the proof.
\epf
\brem\label{better_rate_apost}
From Theorem \ref{qrm_err_est_apost} (i) it follows that if $2<p<b-2$ or if $p>b-2>2$, then the order of the estimate is better than $\mathcal{O}(\d^\frac{1}{2}),$ the best possible rate obtained for the aposteriori case in \cite[Theorem 4.2]{quoc_long_oregan_ngoc_tuan_2022} using the Tikhonov regularization method.

Also, if $4<p<b-2$ or $p>b-2>4,$ the by Theorem \ref{qrm_err_est_apost} (i) it follows that the order of the estimate is better than $\mathcal{O}(\d^\frac{2}{3}).$
\erem
\section{Optimality}\label{sec-optimality}

In this section we analyse the order optimality of the error estimates obtained in the preceding section for the source set $S_{\varrho,p}$. We have already observed in Section \ref{sec_prelim}   that the considered inverse source identification problem is equivalent to solving a linear operator equation 
$$ {\mathbb T}f = h,$$
where ${\mathbb T} : L^2(\O)\to L^2(\O)$ is a compact  self-adjoint bounded linear operator  defined by  (\ref{cpt-op}), that is, 
$${\mathbb T}\f = \sum_{n=1}^\infty \mu_n \<\f, \f_n\>\f_n,\q \f\in L^2(\O)$$
with $\mu_n$  as in (\ref{ev}).

We now formulate some useful representation of the source set $S_{\varrho, p}$ in terms of the operator $\mathbb{T}.$
For this, we observe from  Theorem \ref{non_zero_lower_bound} and  \eqref{bound_for_int_e_psi} that 
\beq\label{lambda-mu} \frac{\lambda_n^2}{\|\psi\|_\infty} \leq \frac{1}{|\mu_n|}\leq \frac{\lambda_n^2}{C}.\eeq 
Thus, in particular we have
\beq\label{upper_lower_bound_for_e_psi}
\frac{\lambda_n^{2p}}{\|\psi\|^p_\infty}\leq \frac{1}{|\mu_n|^p}\leq \frac{\lambda_n^{2p}}{C^p}.
\eeq
Therefore, from \eqref{upper_lower_bound_for_e_psi} it follows that 
\beq\label{equiv_source_cond}
\sum_{n=1}^\infty \lambda_n^{2p}|\<\tilde{f},\f_n\>|^2\leq \varrho^2\q\text{if and only if}\q \sum_{n=1}^\infty \frac{|\<\tilde{f},\f_n\>|^2}{|\mu_n|^p}\leq \tilde{\varrho}^2,
\eeq
where, for a given $\varrho>0$, $\tilde \varrho:= C^{-\frac{p}{2}}\varrho$, and for a given $\tilde \varrho>0$, $\varrho:= \|\psi\|_\infty^{\frac{p}{2}}\tilde \varrho$.
%$\tilde{\varrho}=k\varrho$, and $k\in\{C^{-\frac{p}{2}},\|\psi\|_\infty^{-\frac{p}{2}}\}$ is chosen accordingly. %\hfill {\bf \cred This sentence is not correct}

Next, recall that the operator ${\mathbb T}: L^2(\O)\to L^2(\O)$ defined above 
is a compact self-adjoint operator. Hence, ${\mathbb T}^*{\mathbb T} = {\mathbb T}^2$ is a positive self-adjoint operator with representation
$${\mathbb T}^*{\mathbb T}   \f=\sum_{n=1}^\infty \mu_n^2\<\f,\f_n\>\,\f_n,\q \f\in L^2(\O),$$
and its square-root, conventionally denoted by $|{\mathbb T}|$, has the representation 
$$|{\mathbb T}|   \f=\sum_{n=1}^\infty |\mu_n|\<{\f},\f_n\>\,\f_n,\q \f\in L^2(\O).$$
Consequently, for any $\nu>0$, we have 
$$|{\mathbb T}|^\nu {\f}=\sum_{n=1}^\infty |\mu_n|^\nu\<{\f},\f_n\>\,\f_n,\q \f\in L^2(\O).$$
We may also observe that, for $\nu>0$ if 
$${\mathcal D}_\nu := \{ \f\in L^2(\O): \sum_{n=1}^\infty |\mu_n|^{-2\nu} |\<{\f},\f_n\>|^2 <\infty\},$$
then ${\mathcal D}_\nu$ is the range of $|{\mathbb T}|^\nu$ so that ${\mathcal D}_\nu$ is the domain of the unbounded operator  $|{\mathbb T}|^{-\nu}$, which has the representation 
$$|{\mathbb T}|^{-\nu} \f=\sum_{n=1}^\infty |\mu_n|^{-\nu} \<{\f},\f_n\>\,\f_n,\q \f\in {\mathcal D}_\nu.$$
 
%%Therefore, we have
%%$$(\mathbb{T}^*\mathbb{T})^{-\frac{p}{4}}\,\tilde{f}=\sum_{n=1}^\infty  |\mu_n|^{-\frac{p}{2}}\<\tilde{f},\f_n\>\,\f_n,$$
%%and hence $$\|(\mathbb{T}^*\mathbb{T})^{-\frac{p}{4}}\,\tilde{f}\|^2_{L^2(\O)}=\sum_{n=1}^\infty |\mu_n|^{-p}|\<\tilde{f},\f_n\>|^2.$$
%
%
Therefore, from \eqref{equiv_source_cond}, we have
$$\sum_{n=1}^\infty \lambda_n^{2p}|\<f,\f_n\>|^2\leq \varrho^2\q\text{if and only if}\q \||\mathbb{T}|^{-\frac{p}{2}}f\|^2_{L^2(\O)}\leq \tilde{\varrho}^2$$
whenever $f\in {\mathcal D}_{p/2}$.
Thus, for $r>0$, we have  
\beq\label{source-set} M_{r, p} :=\{g\in L^2(\O):g= |\mathbb{T}|^\frac{p}{2}\tilde{g},\,\|\tilde{g}\|_{L^2(\O)}\leq r \} =  \{g\in {\mathcal D}_{p/2}: \||{\mathbb T}|^{-p/2}g\|_{L^2(\O)} \leq r\}.\eeq 
Now, it is easy to see that $S_{\varrho,p}\subseteq M_{r_1,p}$ and $M_{r_2,p}\subseteq S_{\varrho,p}$, where $r_1=C^{-\frac{p}{2}}\varrho$ and $r_2=\|\psi\|_{\infty}^{-\frac{p}{2}}\,\varrho.$

% Let the {\it worst case error} for a reconstruction method described by an operator $\mathcal{R}:L^2(\O)\to L^2(\O)$ approximating the inverse of $\mathbb{T}$, be denoted by $\triangle(\mathcal{R},\d,M_{\varrho,p})$. Then recall that (cf. \cite{engl_hanke_neubauer, nair_opeq})
%$$\triangle(\mathcal{R},\d,M_{\varrho, p}):=\sup\{\|\mathcal{R}\tilde{h}-\tilde{f}\|: \tilde{f}\in M_{\varrho, p},\,\tilde{h}\in L^2(\O),\,\|\mathbb{T}\tilde{f}-\tilde{h}\|_{L^2(\O)}\leq \d \}.$$
%Also, recall that the {\it modulus of continuity} (cf. \cite{engl_hanke_neubauer, nair_opeq}) denoted by $\o(\d,M_{\varrho, p})$, is defined by
%$$\o(\d,M_{\varrho, p}):=\sup\{\|\tilde{f}_1-\tilde{f}_2\|_{L^2(\O)}: \tilde{f}_1,\tilde{f}_2\in M_{\varrho,p}\,,\, \|\mathbb{T}\tilde{f}_1-\mathbb{T}\tilde{f}_2\|_{L^2(\O)}\leq \d \}.$$ 
%It is well known that the worst case error is bounded below by the modulus of continuity (upto a constant multiple). More precisely, we have 
%\beq\label{uni_lower_bound_for_worst_case_error}
%\triangle(\mathcal{R},\d,M_{\varrho, p})\geq \frac{1}{2}\o(2\d,M_{\varrho, p})\q\text{for all}\q\d>0,
%\eeq
%see for e.g., \cite[Remark 3.12]{engl_hanke_neubauer}. We also observe an easily verifiable fact:
%\beq\label{monotone_prop_mod_cont}
%\o(\d_1,M_{\varrho,p})\leq \o(\d_2,M_{\varrho,p})\q\text{for}\q 0<\d_1\leq \d_2.
%\eeq

We may recall that  (cf. \cite{engl_hanke_neubauer, nair_opeq}), corresponding to $\d>0$,  the {\it worst case error} for a reconstruction method described by an operator $\mathcal{R}:L^2(\O)\to L^2(\O)$ approximating the inverse of $\mathbb{T}$ \w.r.t. a bounded set $\mathcal{M} \subset L^2(\O)$ is defined to be the quantity  
$$\triangle(\mathcal{R},\d,\mathcal{M}):=\sup\{\|\mathcal{R}\tilde{h}-\tilde{f}\|: \tilde{f}\in \mathcal{M},\,\tilde{h}\in L^2(\O),\,\|\mathbb{T}\tilde{f}-\tilde{h}\|\leq \d \}.$$
It is known that (cf. \cite{engl_hanke_neubauer, nair_opeq}) the worst case error is bounded below by the corresponding {\it modulus of continuity} 
$$\o(\d,\mathcal{M}):=\sup\{\|\tilde{f}_1-\tilde{f}_2\|_{L^2(\O)}: \tilde{f}_1,\tilde{f}_2\in \mathcal{M}\,,\, \|\mathbb{T}\tilde{f}_1-\mathbb{T}\tilde{f}_2\|\leq \d \},$$ 
up to a constant multiple, in the sense that %(see for e.g., \cite[Remark 3.12]{engl_hanke_neubauer})
\beq\label{uni_lower_bound_for_worst_case_error}
\triangle(\mathcal{R},\d,\mathcal{M})\geq \frac{1}{2}\o(2\d,\mathcal{M})\q\text{for all}\q\d>0.
\eeq
We also observe some easily verifiable facts:
\beq\label{monotone_prop_mod_cont}
\o(\d_1,\mathcal{M})\leq \o(\d_2,\mathcal{M})\q\text{for}\q 0<\d_1\leq \d_2,
\eeq
and
\beq\label{set_monotone_prop_mod_cont}
\o(\d,\mathcal{M}_1)\leq \o(\d,\mathcal{M}_2)\q\text{for}\q \mathcal{M}_1\subseteq \mathcal{M}_2,\q\d>0.
\eeq
For more details on results related to modulus of continuity for linear operators, we refer to \cite{hofmann_mathe_schieck_2008,nair_opeq}.

We will denote the spectrum of a bounded linear operator $\mathbb{K}$ by ${\rm \bf spec}(\mathbb{K}).$ The following result is motivated from Trong and Hai \cite{trong_hai_2021}. In its proof, we shall make use of the function $\r(\cdot)$ defined by 
$$\r(t):= t^\frac{p+2}{p},\q t>0,$$ 
and the relation 
$\r(t)=t \f^{-1}(t),\, t>0,$
where 
$\f(t):=t^\frac{p}{2}$ for $t>0$. Clearly, $\f$ and $\r$ satisfies the assumption in \cite[Assumption 1.1]{tautenhahn_1998}, that is, $\ds\lim_{t\to 0} \f(t)=0,$ $\f$ is strictly increasing and $\r$ is convex. Note that $|\mathbb{T}|^2\f(|\mathbb{T}|^2)=|\mathbb{T}|^{p+2}$ and $M_{r,p}=\{g\in L^2(\O):g= [\f(|\mathbb{T}|^2)]^\frac{1}{2}\tilde{g},\,\|\tilde{g}\|_{L^2(\O)}\leq r \}$.

\bt\label{minmax_bounds_mod_cont}
Let $\d>0, r>0$ and $M_{r,p}$ be as in (\ref{source-set}). Then the following results hold.
\ben
\i[(i)] If $\frac{\d^2}{r^2}\in {\rm {\bf spec}} \big(|\mathbb{T}|^{p+2}\big),$ then $\o(\d,M_{r,p})= r^\frac{2}{p+2}\d^\frac{p}{p+2}.$
\i[(ii)] If $\frac{\d^2}{r^2}\notin {\rm {\bf spec}} \big(|\mathbb{T}|^{p+2}\big)$, then for $\d\leq r \sup\limits_{n}|\mu_n|^\frac{p+2}{2}$ either
$$ \inf_{n}\frac{|\mu_n|^{\frac{p}{2}}}{|\mu_{n+1}|^{\frac{p}{2}}}\,r^\frac{2}{p+2}\,\d^\frac{p}{p+2}\leq \o(\d,M_{r,p})\leq \sup_{n}\frac{|\mu_{n+1}|^\frac{p}{2}}{|\mu_n|^\frac{p}{2}}\,r^\frac{2}{p+2}\,\d^\frac{p}{p+2},$$
or
$$\inf_{n}\frac{|\mu_{n+1}|^{\frac{p}{2}}}{|\mu_n|^{\frac{p}{2}}}\,r^\frac{2}{p+2}\,\d^\frac{p}{p+2}\leq \o(\d,M_{r,p})\leq \sup_{n}\frac{|\mu_n|^\frac{p}{2}}{|\mu_{n+1}|^\frac{p}{2}}\,r^\frac{2}{p+2}\,\d^\frac{p}{p+2}.$$
\een
\et

\bpf
\noi {\bf Case 1.}\q $\frac{\d^2}{r^2}\in {\rm {\bf spec}} \big(|\mathbb{T}|^{p+2}\big).$

In this case, by the result of Tautenhahn \cite[Theorem 2.1]{tautenhahn_1998}, it follows that $$\o(\d,M_{r,p})=r \sqrt{\r^{-1}\left(\frac{\d^2}{r^2}\right)}=r\left(\frac{\d}{r}\right)^\frac{p}{p+2}=r^\frac{2}{p+2}\,\d^\frac{p}{p+2}.$$
{\bf Case 2.}\q $\frac{\d^2}{r^2}\notin {\rm {\bf spec}} \big(|\mathbb{T}|^{p+2}\big).$

Recall that the set of eigenvalues of the compact self-adjoint operator $|\mathbb{T}|^{p+2}$ is precisely the set $\Big\{|\mu_n|^{p+2}:n\in \N\Big\},$
and that $|\mu_n|\to 0$ as $n\to \infty.$ Thus,  there exists $n_0\in \N$ such that $|\mu_{n_0}|=\sup\limits_{n}\Big\{|\mu_n|:n\in \N\Big\}.$
Therefore, for $\frac{\d^2}{r^2}\notin {\rm {\bf spec}} \big(|\mathbb{T}|^{p+2}\big)$ and for $0<\d\leq r \sup\limits_{n}|\mu_n|^\frac{p+2}{2}$, there are only two possibilities given below as sub-cases.

\noi {\bf \underline{Sub-case 2(i).}}\q There exists $n\in \N$ such that 
$$|\mu_n|^{p+2}\leq \frac{\d^2}{r^2}\leq |\mu_{n+1}|^{p+2},$$
that is, 
\beq\label{inter_ineq_optimal}
r\, |\mu_n|^{\frac{p+2}{2}}\leq \d\leq r \,|\mu_{n+1}|^\frac{p+2}{2}.
\eeq
Therefore, by \eqref{monotone_prop_mod_cont} and \eqref{inter_ineq_optimal}, we have
\beq\label{omega_lower_bound_case21}
\o(\d,M_{r,p})\geq \o(r\,|\mu_n|^{\frac{p+2}{2}},M_{r,p})\underbrace{=}_{{\bf Case 1}} r^\frac{2}{p+2}\left(r\,|\mu_n|^{\frac{p+2}{2}}\right)^\frac{p}{p+2}.\eeq
From \eqref{inter_ineq_optimal} we have
$$\frac{r}{\d}\,|\mu_{n+1}|^\frac{p+2}{2}\geq 1$$ and hence by \eqref{omega_lower_bound_case21}
$$\o(\d,M_{r,p})\geq r^\frac{2}{p+2}\,\d^\frac{p}{p+2}\,\frac{|\mu_n|^{\frac{p}{2}}}{|\mu_{n+1}|^{\frac{p}{2}}}.$$
Again from \eqref{monotone_prop_mod_cont} and \eqref{inter_ineq_optimal}, we have
\beq\label{omega_upper_bound_case21}
\o(\d,M_{r,p})\leq \o(r\,|\mu_{n+1}|^\frac{p+2}{2}, M_{r,p})\underbrace{=}_{{\bf Case 1}} r^\frac{2}{p+2}\,\left(r\,|\mu_{n+1}|^\frac{p+2}{2}\right)^\frac{p}{p+2}.
\eeq
Now from \eqref{inter_ineq_optimal}, we have
$$\frac{r}{\d}\,|\mu_n|^\frac{p+2}{2}\leq 1$$
and hence by \eqref{omega_upper_bound_case21}, we have
$$\o(\d,M_{r,p})\leq r^\frac{2}{p+2}\,\d^\frac{p}{p+2}\,\frac{|\mu_{n+1}|^\frac{p}{2}}{|\mu_n|^\frac{p}{2}}.$$
Thus, we have
$$ \inf_{n}\frac{|\mu_n|^{\frac{p}{2}}}{|\mu_{n+1}|^{\frac{p}{2}}}\,r^\frac{2}{p+2}\,\d^\frac{p}{p+2}\leq \o(\d,M_{r,p})\leq \sup_{n}\frac{|\mu_{n+1}|^\frac{p}{2}}{|\mu_n|^\frac{p}{2}}\,r^\frac{2}{p+2}\,\d^\frac{p}{p+2}.$$
{\bf \underline{Sub-case 2(ii).}}\q There exists $n\in \N$ such that 
$$|\mu_{n+1}|^{p+2}\leq \frac{\d^2}{r^2}\leq |\mu_n|^{p+2}.$$
Then using similar arguments as done for {\bf Sub-case 2(i)}, we obtain
$$\inf_{n}\frac{|\mu_{n+1}|^{\frac{p}{2}}}{|\mu_n|^{\frac{p}{2}}}\,r^\frac{2}{p+2}\,\d^\frac{p}{p+2}\leq \o(\d,M_{r,p})\leq \sup_{n}\frac{|\mu_n|^\frac{p}{2}}{|\mu_{n+1}|^\frac{p}{2}}\,r^\frac{2}{p+2}\,\d^\frac{p}{p+2}.$$
\epf
\brem\label{delta_0}
From \eqref{bound_for_int_e_psi} it follows that $\sup\limits_n|\mu_n| \leq \frac{\|\psi\|_\infty}{\lambda_1^2}.$ Therefore, Theorem \ref{minmax_bounds_mod_cont} (ii) holds for all $0<\d\leq \d_0,$ where $\d_0:=r\,\sup\limits_n|\mu_n|^\frac{p+2}{2} < \infty.$%{\cb r}\left(\frac{\|\psi\|_\infty}{\lambda_1^2}\right)^\frac{p+2}{2}.$
\erem
\brem
We now show that the $infimums$ and the $supremums$ in Theorem \ref{minmax_bounds_mod_cont}(ii) are indeed non-zero and finite, respectively. In order to see this, first we recall a  result from Courant and Hilbert \cite{courant_hilbert_book}, which says that there exist constants $e_1, e_2>0$ such that $$e_1 n^\frac{2}{d}\leq \lambda_n\leq e_2 n^\frac{2}{d}\q\text{for all}\q n\in \N.$$
From this, we obtain 
$$ \frac{e_1}{e_2} \Big(\frac{1}{2}\Big)^{\frac{2}{d}} \leq \frac{e_1}{e_2} \Big(\frac{n}{n+1}\Big)^{\frac{2}{d}}  \leq \frac{\l_n}{\l_{n+1}} \leq \frac{e_2}{e_1} \Big(\frac{n}{n+1}\Big)^{\frac{2}{d}} \leq \frac{e_2}{e_1}.$$
This together with the inequalities in \eqref{lambda-mu} imply the required assertion. 
%
%
%
%This leads to $$\lambda_n\leq \lambda_{n+1}\leq \frac{e_2\lambda_n}{e_1}\left(1+\frac{1}{n}\right)^\frac{2}{d}\leq \lambda_n\frac{2^{2/d}\,e_2 }{e_1}\q\text{for all}\q n\in\N.$$
%Now, 
%\beqarray
%\inf_{n}\frac{|\mu_n|}{|\mu_{n+1}|}&\underbrace{\geq}_{\eqref{bound_for_int_e_psi}} &\frac{1}{\|\psi\|_\infty}\inf_{n}\, \lambda_{n+1}^2\,|\mu_n|
%\underbrace{\geq}_{{\rm Theorem} \ref{non_zero_lower_bound}} \frac{C}{\|\psi\|_\infty}\,\inf_{n}\,\frac{\lambda_{n+1}^2}{\lambda_n^2}\geq \frac{C}{\|\psi\|_\infty}
%\eeqarray
%and 
%\beqarray
%\inf_{n}\frac{|\mu_{n+1}|}{|\mu_n|} &\underbrace{\geq}_{\eqref{bound_for_int_e_psi}} & \frac{1}{\|\psi\|_{\infty}}\,\inf_{n}\,\lambda_n^2\,|\mu_{n+1}|
%\underbrace{\geq}_{{\rm Theorem} \ref{non_zero_lower_bound}}  \frac{C}{\|\psi\|_\infty}\,\inf_{n}\,\frac{\lambda_n^2}{\lambda_{n+1}^2} \geq \frac{C}{\|\psi\|_\infty}\,\left(\frac{e_1}{2^{2/d}\,e_2}\right)^2.
%\eeqarray
%Similarly, we obtain
%$$\sup_{n}\frac{|\mu_{n+1}|}{|\mu_n|}\leq \frac{\|\psi\|_\infty}{C}\q\h{and}\q 
%\sup_{n}\frac{|\mu_n|}{|\mu_{n+1}|}\leq \frac{\|\psi\|_\infty}{C}\,\left( \frac{2^{2/d}\,e_2}{e_1}\right)^2$$
\erem

\brem{\rm (\textbf{Order optimality})}\label{order_opt}
Let $\triangle(QRM, \d,S_{\varrho,p})$ denote the worst case error for the quasi-reversibilty method. Recall that $M_{r_2,p}\subseteq S_{\varrho,p},$ where $r_2=\|\psi\|_\infty^{-\frac{p}{2}}\varrho.$ Let $0<\d\leq \d_0,$ where $\d_0:= \frac{r_2}{2}\, \sup\limits_{n}|\mu_n|^\frac{p+2}{2}$. 

Then by Theorem \ref{qrm_err_est_apriori} (i) for $0<p<b$ in the {\bf apriori case}, and for any reconstruction method $\mathcal{R}$, we have
\beqarray
\triangle(QRM,\d,S_{\varrho,p})&\leq & C_1\,\varrho^\frac{2}{p+2}\,\d^\frac{p}{p+2}\, =\frac{C_1}{C_\o}\, C_\o\,\varrho^\frac{2}{p+2}\,\d^\frac{p}{p+2}\\
& \underbrace{\leq}_{{\rm Theorem}\ref{minmax_bounds_mod_cont}}& \frac{C_1}{C_\o}\,\o(2\d,M_{r_2,p})\underbrace{\leq}_{\eqref{set_monotone_prop_mod_cont}} \frac{C_1}{C_\o}\,\o(2\d,S_{\varrho,p})\underbrace{\leq}_{\eqref{uni_lower_bound_for_worst_case_error}} \frac{2C_1}{C_\o}\, \triangle(\mathcal{R},\d,S_{\varrho,p}),
\eeqarray
where $C_\o$ is accordingly chosen from the set $$\Big\{2^\frac{p}{p+2}\|\psi\|_\infty^{-\frac{p}{p+2}},\q \inf_{n}\frac{|\mu_n|^{\frac{p}{2}}}{|\mu_{n+1}|^{\frac{p}{2}}}\,2^\frac{p}{p+2}\|\psi\|_\infty^{-\frac{p}{p+2}},\q \inf_{n}\frac{|\mu_{n+1}|^{\frac{p}{2}}}{|\mu_n|^{\frac{p}{2}}}\,2^\frac{p}{p+2}\|\psi\|_\infty^{-\frac{p}{p+2}}\Big\}.$$
Similarly, by Theorem \ref{qrm_err_est_apost} (i) for $b>2$ and $0<p<b-2$ in the {\bf aposteriori case}, and for any reconstruction method $\mathcal{R},$ we have
$$\triangle(QRM,\d,S_{\varrho,p})\leq \frac{2C_6}{C_\o}\, \triangle(\mathcal{R},\d, S_{\varrho,p}).$$

Therefore, in the apriori case for $0<p<b$ and in the aposteriori case for $b>2$ and $0<p<b-2$ the error estimates obtained in Theorem \ref{qrm_err_est_apriori} and Theorem \ref{qrm_err_est_apost}, respectively, are of optimal order for the source set $S_{\varrho,p}$.
\erem
\brem\label{incomp_arg_duc_hao}
In \cite[Theorem 2]{duc_thang_thanh_2023} and \cite[Remark 2]{hao_liu_duc_thang_2019} the authors have claimed the order optimality of their estimates in respective context by appealing to the result of Tautenhahn \cite[Theorem 2.1]{tautenhahn_1998}, which relies on the fact that $\frac{\d^2}{\varrho^2}$ is in the spectrum of the specified operators under consideration. Since this condition need not be satisfied for a  given $\d>0$, however small it may be, the arguments for the order optimality in \cite[Theorem 2]{duc_thang_thanh_2023} and \cite[Remark 2]{hao_liu_duc_thang_2019} is obviously incomplete.
%, because of not considering the case when $\frac{\d^2}{\varrho^2}$ does not belongs to those spectral set. 
Our arguments for Theorem \ref{minmax_bounds_mod_cont} (ii) and Corollary \ref{order_opt} can be modified appropriately in those settings to fill the gaps in their arguments.
\erem
\section{conclusion}
We have considered an inverse source identification problem from final time observation associated with a bi-parabolic system, which is known to be an ill-posed problem. In order to obtain stable approximations of the sought source, we have employed the quasi-reversibilty (QRM) method, which has not been explored yet. We obtained H{\"o}lder type error estimates for both the apriori and aposteriori parameter choice strategies. Our rates seems to exceed the rates obtained in earlier works by using Tikhonov regularization and quasi-boundary value method. Moreover, we have also shown the order optimality of the obtained rates for some cases. The arguments used for proving order optimality in this paper, can be used to fill the gap in the arguments of order optimality claimed in some earlier works in different context by different authors.

Moreover, this paper seems to be the first one that has broaden the applicability of the considered problem by enlarging the admissible class of source functions. This is done so by allowing the time dependent component of the source function to change sign. Note that in the context of the considered inverse problem, all previously known works have assumed some fixed sign property of the time dependent component by assuming some bounded below conditions.

\vspace{1cm}
\noi{\bf Acknowledgments.} The first author, Subhankar Mondal,  is supported by the postdoctoral fellowship of TIFR Centre for Applicable Mathematics, Bangalore, and the second author, M. Thamban Nair,  gratefully acknowledges the
support received from BITS Pilani, K.K. Birla Goa Campus, where he is a Visiting
Professor.

\end{document}